\documentclass[12pt]{amsart}
\usepackage[top=1.08in, left=1.28in, bottom=1.08in, right=1.28in]{geometry}
\usepackage{amsmath,amsfonts,amsbsy,amsgen,amscd,mathrsfs,amssymb,amsthm,mathtools,bbm,enumerate}
\usepackage[colorlinks=true,citecolor=teal!80!black,linkcolor=blue]{hyperref}
\usepackage[foot]{amsaddr}

\usepackage{hhline}
\usepackage{fouriernc}

\usepackage{pgfplots}
\usepackage{tikz-cd}
\usetikzlibrary{arrows,automata}

\numberwithin{equation}{section}
\newtheorem{theorem}{Theorem}[section]
\newtheorem{corollary}[theorem]{Corollary}

\newtheorem{lemma}[theorem]{Lemma}
\newtheorem{proposition}[theorem]{Proposition}
\theoremstyle{definition}

\newtheorem{definition}[theorem]{Definition}
\newtheorem{example}[theorem]{Example}
\newtheorem{heuristic}[theorem]{Heuristic}

\newtheorem{openproblem}[theorem]{Open Problem}
\newtheorem{question}[theorem]{Question}
\newtheorem{remark}[theorem]{Remark}

\makeatletter

\newcommand\al{\alpha}
\newcommand\be{\beta}
\newcommand\dd{\mathrm d}
\newcommand\De{\Delta}

\newcommand\deq{\stackrel{\mathrm{distr.}}{=}}
\newcommand\eps{\varepsilon}
\newcommand\ga{\gamma}

\newcommand\ka{\kappa}
\newcommand\La{\Lambda}
\newcommand\la{\lambda}

\newcommand\Si{\Sigma}
\newcommand\si{\sigma}

\renewcommand\d{~\mathrm d}
\renewcommand\phi{\varphi}
\renewcommand\rho{\varrho}

\renewcommand\hat{\widehat}

\newcommand\bs{\boldsymbol}
\newcommand\mbb{\mathbb}
\newcommand\mbf{\mathbf}
\newcommand\mc{\mathcal}
\newcommand\mf{\mathfrak}
\newcommand\mr{\mathrm}

\newcommand\msf{\mathsf}

\begin{document}

\title[Recovering Parameters from Edge Fluctuations]{Recovering Parameters from Edge Fluctuations: Beta-Ensembles and Critically-Spiked Models}
\author{Pierre Yves Gaudreau Lamarre}
\address{Department of Mathematics,
Syracuse University,
Syracuse, NY 13244}
\email{pgaudrea@syr.edu}
\maketitle

\begin{abstract}
Let $\Lambda=\{\Lambda_0,\Lambda_1,\Lambda_2,\ldots\}$ be the point process that describes the edge scaling limit
of either (i) "regular" beta-ensembles with inverse temperature $\beta>0$, or
(ii) the top eigenvalues of Wishart or Gaussian invariant random matrices perturbed by $r_0\geq1$ critical spikes.
In other words, $\Lambda$ is the eigenvalue point process of one of the scalar or multivariate stochastic Airy operators.
We prove that a single observation of $\Lambda$ suffices to recover (almost surely) either (i) $\beta$ in the case of beta-ensembles,
or (ii) $r_0$ in the case of critically-spiked models. Our proof relies on the recently-developed semigroup theory
for the multivariate stochastic Airy operators.

Going beyond these parameter-recovery applications, our results also (iii) refine our understanding
of the rigidity properties of $\Lambda$,
and (iv) shed new light on the equality (in distribution) of stochastic Airy spectra with different dimensions and
the same Robin boundary conditions.
\end{abstract}

\section{Introduction}

\subsection{Main Result}

Let $\Theta$ be the set of parameters of the form $\theta=(r,\be,w)$, where $r\in\mbb N$ is any positive integer and
$\be$ and $w$ are defined as follows:
\begin{enumerate}[(1)]
\item If $r=1$, then $\be>0$ and $w\in(-\infty,\infty]$, and
\item if $r>1$, then $\be\in\{1,2,4\}$ and $w=(w_1,\ldots,w_r)\in(-\infty,\infty]^r$.
\end{enumerate}
Given $\theta\in\Theta$, let $\mc H_{\theta}$ be the corresponding stochastic Airy operator (SAO). That is:
\begin{enumerate}[(1)]
\item If $r=1$, then we let $W:[0,\infty)\to\mbb R$ be a standard Brownian motion, and let
\begin{align}
\label{Equation: Scalar SAO}
\mc H_\theta=-\tfrac{\mr d^2}{\mr d x^2}+x+\tfrac{2}{\be^{1/2}}W'(x),
\end{align}
which acts on functions $f\in L^2\big([0,\infty),\mbb R\big)$ subject to the boundary condition
\[\begin{cases}
f(0)=0&\text{if }w=\infty\text{ (Dirichlet)},\\
f'(0)=wf(0)&\text{if }w\in\mbb R\text{ (Robin)}.
\end{cases}\]
\item Let $\mbb F_1=\mbb R$,
$\mbb F_2=\mbb C$, and $\mbb F_4=\mbb H$.
If $r>1$ and $\be\in\{1,2,4\}$, then
let $W_\be$ be a standard $r$-dimensional matrix Brownian motion on $\mbb F_\be$
(i.e., with GOE, GUE, or GSE increments depending on whether $\be=1,2,4$; see \cite[Page 2730]{BloemendalVirag2}), and
define
\begin{align}
\label{Equation: Multivariate SAO}
\mc H_\theta=-\tfrac{\mr d^2}{\mr d x^2}+rx+\sqrt 2W_\be'(x).
\end{align}
This acts on $r$-dimensional vector-valued functions $f=\big(f(i,\cdot):1\leq i\leq r\big)$ whose components
$f(i,\cdot)\in L^2\big([0,\infty),\mbb F_\be\big)$ are subject to the boundary conditions
\begin{align}
\label{Equation: wi}
\begin{cases}
f(i,0)=0&\text{if }w_i=\infty\text{ (Dirichlet)},\\
f'(i,0)=w_if(i,0)&\text{if }w_i\in\mbb R\text{ (Robin)},
\end{cases}\qquad 1\leq i\leq r.
\end{align}
\end{enumerate}

The SAOs were introduced by Edelman-Sutton \cite{EdelmanSutton},
Ram\'irez-Rider-Vir\'ag \cite{RamirezRiderVirag}, and Bloemendal-Vir\'ag \cite{BloemendalVirag1,BloemendalVirag2}
to describe the edge fluctuations of various point processes of high interest in mathematical physics
and statistics. That is, depending on the choice of $\theta$, the
eigenvalues of $\mc H_\theta$ either describe the edge fluctuations of a wide class of 
beta-ensembles with so-called "regular" external potential
\cite{BekermanFigalliGuionnet,Bekerman,BourgadeErdosYau,DeiftGioev,KrishnapurRiderVirag,PasturShcherbina,RamirezRiderVirag,Shcherbina},
or the soft-edge fluctuations of critically-spiked Wishart or Gaussian invariant ensembles with real, complex, or quaternion entries \cite{BaikBenArousPeche,BloemendalVirag1,BloemendalVirag2,Mo,Peche,Wang}.

In this paper, we are interested in the following question regarding the SAOs:

\begin{question}
Let $\theta\in\Theta$ be fixed, and let
$\La^\theta=\big\{\La^\theta_0,\La^\theta_1,\La^\theta_2,\ldots\big\}$
denote the eigenvalues of $\mc H_\theta$ arranged in increasing order.
What information about $\theta$ can be recovered almost-surely from
a single realization of $\La^\theta$?
\end{question}

Given $(r,\be,w)\in\Theta$, define the parameter
\begin{align}
\label{Equation: r0}
r_0=r_0(w)=\sum_{i=1}^r\mbf 1_{\{w_i<\infty\}};
\end{align}
if $r=1$, then this reduces to $r_0=\mbf 1_{\{w<\infty\}}$.
In words, $r_0$ is the number of Robin boundary conditions imposed on $\mc H_\theta$'s domain
(hence $r-r_0$ is the number of Dirichlet conditions).
Our main result asserts that if either one of $\be$ or $r_0$ is known,
then the other parameter can be recovered almost surely from a single realization of $\La^\theta$:

\begin{theorem}
\label{Theorem: Main Informal}
There exists a deterministic function $\mc T$
(which can be written as an explicit limit; see \eqref{Equation: Main} for the details) such that for every $\theta\in\Theta$, one has
\begin{align}
\label{Equation: Main Informal}
\mc T(\La^\theta)=r_0+1/\be\qquad\text{almost surely.}
\end{align}
\end{theorem}

\subsection{Applications}

Our main motivation for proving Theorem \ref{Theorem: Main Informal} consists of two parameter-recovery results:
\begin{enumerate}[(i)]
\item Corollary \ref{Corollary: Beta Ensembles}: The temperature of any sequence of "regular" beta-ensembles can be recovered from a single
observation of their edge fluctuations.
\item Corollary \ref{Corollary: Spiked}: The number of signals "close enough" to the critical threshold in sequences of
spiked Wishart and Gaussian invariant models can be recovered
from a single observation of their top edge fluctuations.
\end{enumerate}
Part of the significance of these results is that, to the best of our knowledge, all previously-known techniques to recover
these parameters use the asymptotics of the full configuration of the point process or random matrix eigenvalues in question. In sharp contrast to this, Corollaries
\ref{Corollary: Beta Ensembles} and \ref{Corollary: Spiked} show that the same can be done using only the vanishing proportion
of particles/eigenvalues that contribute to edge scaling limits.
See Sections \ref{Section: Beta-Ensembles} and \ref{Section: Spiked Models} for more details and
background, as well as open problems motivated by these results
(Open Problems \ref{Open Problem: Beta Ensembles} and \ref{Open Problem: Spiked}).

In addition to the above, we present two applications of a more theoretical nature:
\begin{enumerate}[(i)]
\setcounter{enumi}{2}
\item In Corollaries \ref{Corollary: Rigidity Explanation} and \ref{Corollary: Parameter Rigidity}, we explain how Theorem \ref{Theorem: Main Informal} refines our understanding of the rigidity properties of $\La^\theta$ (notably number rigidity).
\item In Corollary \ref{Corollary: Cancellations}, we explain how our proof technique
sheds new light on the mysterious invariance of SAO eigenvalues with respect to adding
components with Dirichlet boundary conditions.
\end{enumerate}
We refer once again to Sections \ref{Section: Rigidity} and \ref{Section: Cancellations} for more details, background, and open problems motivated by these applications.

\subsection{Proof Technique}

Our approach to prove Theorem \ref{Theorem: Main Informal} is inspired by the theory
of inverse problems in spectral geometry. That is, the observation
that an impressive amount of information
about a differential operator's domain can be inferred from its eigenvalues alone.
One of the most well-known manifestation of this phenomenon is the Weyl law
\cite{Weyl} and its various extensions: If $\big\{\la_k(M,b):k\geq0\big\}$ are the eigenvalues of the Laplacian
on some compact manifold $M$ subject to some boundary condition $b$ (if $\partial M\neq\varnothing$),
then $M$'s dimension, $M$'s volume, $M$'s boundary area, and the boundary condition $b$ (at least Dirichlet versus Neumann)
can in many cases be determined by the asymptotics of the eigenvalue counting function
\[\la\mapsto\big|\big\{k\geq0:\la_k(M,b)\leq\la\big\}\big|.\]

In this context, Question \ref{Question: Beta Ensemble Recovery} can be viewed as
the inverse spectral problem for the SAOs, as
the components of $\theta=(r,\be,w)$ characterize
$\mc H_\theta$'s domain and noise. Thus,
the proof of Theorem \ref{Theorem: Main Informal} relies on the observation that
$r_0+1/\be$ can be recovered from precise asymptotics of the $\La_k^\theta$'s as $k\to\infty$.
That being said, instead of recovering $r_0+1/\be$ using the asymptotics of the eigenvalue
counting function, we do so
via the exponential traces $\mr{Tr}\big[\mr e^{-t\mc H_\theta/2}\big]=\sum_{k=0}^\infty\mr e^{-t\La^\theta_k/2}$ for small $t>0$
(which is also widely used in spectral geometry; see, e.g., \cite[Chapter 3]{Gilkey} and references therein).
More specifically, our proof of Theorem \ref{Theorem: Main Informal} relies on showing that
\begin{align}
\label{Equation: Informal Trace Asymptotics}
\mr{Tr}\big[\mr e^{-t\mc H_\theta/2}\big]=\sqrt{\frac2\pi}t^{-3/2}+\frac12\left(r_0+\frac1\be\right)-\frac14+\mf o_t\qquad\text{for all }t\in(0,1],
\end{align}
where the random remainder terms $\mf o_t$ are such that
\begin{align}
\label{Equation: Informal Trace Asymptotics Condition 1}
\mbf E[\mf o_t]=o(1)\qquad\text{ as }t\to0^+,
\end{align}
and such that there exists a constant $C>0$ satisfying
\begin{align}
\label{Equation: Informal Trace Asymptotics Condition 2}
\mbf{Cov}[\mf o_s,\mf o_t]\leq C\left(\frac{\min\{s,t\}}{\max\{s,t\}}\right)^{1/4}
\qquad\text{ for every }s,t\in(0,1].
\end{align}
In particular, the fact that any sequence $\mf o_{t_1},\mf o_{t_2},\mf o_{t_3},\ldots$ such that $t_n\to0^+$ has a bounded
variance and vanishing covariance
means that $r_0+1/\be$ can be recovered from \eqref{Equation: Informal Trace Asymptotics} by
the law of large numbers, that is, by averaging the traces $\mr{Tr}\big[\mr e^{-t\mc H_\theta/2}\big]$
over a sequence of $t$'s that vanishes fast enough; see \eqref{Equation: Main} for details.

The main technical input in the proof of \eqref{Equation: Informal Trace Asymptotics}
consists of the recently-developed
semigroup theory for the {\it multivariate} SAOs \cite{RigidityMSAO}, i.e., the case $r>1$ (which is essential for our application
of Theorem \ref{Theorem: Main Informal} to the spiked models in Corollary \ref{Corollary: Spiked}). This development extended prior results on
the semigroup theory of the scalar SAOs \cite{GaudreauLamarreShkolnikov,GorinShkolnikov}, i.e., the case $r=1$.

\subsection{Questions of Optimality}

Since Theorem \ref{Theorem: Main Informal} does not allow to recover the parameter
$\theta$ completely, it is natural to ask if our result is optimal. For instance, one might ask
if, in addition to $r_0$, it is possible to recover information about the magnitudes of the components
$w_i$ in \eqref{Equation: wi} that are finite. That said, it can be shown that this cannot be achieved
with the method used in this paper.

More specifically, it is known that for some cases of $\theta$ (notably, $\theta=(1,2,\infty)$; see \cite[Proposition 2.27]{RigidityFK}),
the random remainder term $\mf o_t$ in \eqref{Equation: Informal Trace Asymptotics} satisfies
\[\liminf_{t\to0^+}\,\mbf{Var}[\mf o_t]>0.\]
Given these nontrivial fluctuations, the deterministic terms in \eqref{Equation: Informal Trace Asymptotics} represent
the full extent of the information that can be inferred with 100\% accuracy from the small $t$ asymptotics
of $\mr{Tr}\big[\mr e^{-t\mc H_\theta/2}\big]$. In particular, if it is in fact possible to recover more information
about $\theta$ from $\La^\theta$, then we expect that a different approach than the one used in this
paper would be required.

\subsection{Organization}

The remainder of this paper is organized as follows:
In Section \ref{Section: Applications}, we provide more details on
the applications and open problems coming from Theorem \ref{Theorem: Main Informal}.
Then, in Section \ref{Section: Proof}, we prove our main result.

\subsection{Acknowledgements}

The author thanks Promit Ghosal and Sumit Mukherjee for
insightful conversations and references regarding the consistent
estimation of inverse temperature in spin glasses and more general Gibbs point processes.

\section{Applications and Open Problems}
\label{Section: Applications}

\subsection{Temperature Recovery in Beta-Ensembles}
\label{Section: Beta-Ensembles}

\subsubsection{Background}

In statistical physics, the beta-ensembles model particle systems on $\mbb R$ subjected to
an external potential energy and pairwise logarithmic repulsion.
More specifically, let $V:\mbb R\to\mbb R$ be a fixed function with sufficient growth, and for every
$n\in\mbb N$, define the $n$-particle Hamiltonian
\begin{align}
\label{Equation: Hamiltonian}
H_n(x)=\frac12\sum_{i=1}^nV(x_i)-\frac1n\sum_{1\leq i<j\leq n}\log|x_j-x_i|,\qquad x=(x_1,\ldots,x_n)\in\mbb R^n.
\end{align}
For every $\be>0$ and $n\in\mbb N$, we define the $n$-dimensional beta-ensemble with external potential $V$ and inverse temperature $\be$
as the point process $\mc B^{\be,n}=\big\{\mc B^{\be,n}_1,\ldots,\mc B^{\be,n}_n\big\}$
sampled according to the Gibbs measure with density
\begin{align}
\label{Equation: Beta Ensemble Density}
\tfrac1{Z_{\be,n}}\mr e^{-\be nH_n(x)}\d x,\qquad x\in\mbb R^n.
\end{align}
We adopt the convention that the points in $\mc B^{\be,n}$ are in increasing order.
The constant $Z_{\be,n}>0$, called the partition function, ensures that \eqref{Equation: Beta Ensemble Density} integrates to one.

The question of whether the temperature of infinite Gibbs point processes
can be recovered from a single sample
has attracted some attention in the literature;
see, e.g.,
\cite[Section 5.4]{Dereudre} and \cite{DereudreLavancier} and references therein.
In the case of a finite configuration, such as \eqref{Equation: Beta Ensemble Density}, this is of course not possible due to the
mutual absolute continuity of the particles' joint densities with different $\be$'s. However, the problem becomes interesting
if one asks about recovering $\be$ as the number of particles grows to infinity.
For instance, several prior works have studied the consistency of maximum likelihood and
pseudolikelihood estimators (MLEs and MPLEs) of temperature in spin glasses in the large-dimensional limit;
see, e.g., \cite{BhattacharyaMukherjee,Chatterjee,CometsGidas,GhosalMukherjee} and references therein.
In view of extending this to beta-ensembles, one may then ask:

\begin{question}
\label{Question: Beta Ensemble Recovery}
Let $V:\mbb R\to\mbb R$ in \eqref{Equation: Hamiltonian} be fixed.
Are there asymptotics (as $n\to\infty$)  of the particles $\mc B^{\be,n}_1\leq\cdots\leq\mc B^{\be,n}_n$ from which
$\be$ can be recovered almost surely? If so, what is the smallest number of particles needed to recover $\be$?
\end{question}

To the best of our knowledge, this has not been studied in the literature.
That being said, one of the main hurdles in applying MLEs to Gibbs point processes
consists of computing the partition function (e.g., \cite[Section 1]{Younes}). Therefore,
given that $Z_{\be,n}$ is very well understood for beta-ensembles---see
\cite{BorotGuionnet1,BorotGuionnet2} and references therein---we expect that
a similar strategy could be used to solve Question \ref{Question: Beta Ensemble Recovery}.
To make this more concrete, consider the following heuristic (which
is informal as the asymptotics therein were calculated using Wolfram Mathematica;
see Appendix \ref{Section: beta-ensembles heuristic}):

\begin{heuristic}[Informal]
\label{Heuristic: beta-ensembles}
Consider the particle configuration's total energy
\begin{align}
\label{Equation: Beta Ensembles Energy}
H_n(\mc B^{\be,n})=\frac12\sum_{i=1}^nV(\mc B^{\be,n}_i)-\frac1n\sum_{1\leq i<j\leq n}\log\big|\mc B^{\be,n}_j-\mc B^{\be,n}_i\big|.
\end{align}
If $V(x)=x^2/2$, then
\begin{align}
\label{Equation: Energy Limit in Probability}
\lim_{n\to\infty}-4\big(H_n(\mc B^{\be,n})-\tfrac38n+\tfrac12\log n\big)-1=\log(\be^2/4)-2\tfrac{\Gamma'(1+\be/2)}{\Gamma(1+\be/2)}\qquad\text{in probability},
\end{align}
where, $\Gamma(z)=\int_0^\infty t^{z-1}\mr e^{-t}\d t$ denotes the Gamma function.
Given that $\be\mapsto\log(\be^2/4)-2\tfrac{\Gamma'(1+\be/2)}{\Gamma(1+\be/2)}$ is invertible for $\be\in(0,\infty)$
(see Appendix \ref{Section: beta-ensembles heuristic}),
we conclude that $\be$ can, in principle, be recovered from the limit \eqref{Equation: Energy Limit in Probability}.
\end{heuristic}

As argued in Appendix \ref{Section: beta-ensembles heuristic}, the limit
\eqref{Equation: Energy Limit in Probability} relies on an asymptotic analysis of the partition function.
Thus, we expect that \eqref{Equation: Energy Limit in Probability} could be extended to more general $V$'s using the asymptotic expansions of
general beta-ensemble partition functions in \cite[(1-3)]{BorotGuionnet1} and \cite[(1.14)]{BorotGuionnet2}.

\subsubsection{New Result}

A common feature of the statistics discussed above (i.e., MLEs, MPLEs, or the total energy in \eqref{Equation: Beta Ensembles Energy}) is that they use
{\it every point} in the configuration $\mc B^{\be,n}_1\leq\cdots\leq\mc B^{\be,n}_n$ as an input. In sharp contrast to this, Theorem \ref{Theorem: Main Informal}
implies that the temperature of many beta-ensembles can be recovered from a much more parsimonious asymptotic,
which only depends on the vanishing proportion of particles that contribute to edge scaling limits (as the vast majority of particles instead contribute to the bulk scaling limits;
e.g., \cite{ValkoVirag,ValkoVirag2}).
This is based on the edge universality of beta-ensembles, part of which can be informally summarized as follows:

\begin{theorem}[Informal]
\label{Theorem: Edge Universality of Beta-Ensembles}
Let $V$ be a regular one-cut potential. That is, $V$ has an equilibrium measure $\mu_V$
supported on a compact interval $[\mf a,\mf b]$, and $\mu_V$'s density vanishes like a square root at
$\mf a$ and $\mf b$. Under additional technical conditions on $V$,
for any edge $\mf e\in\{\mf a,\mf b\}$, there exists scaling constants $\mf s_n=O(n^{2/3})$
such that
\begin{align}
\label{Equation: Beta-Ensembles Scaling Limit}
\lim_{n\to\infty}\mf s_n\big(\mf e-\mc B^{\be,n}\big)=\La^{(1,\be,\infty)}\qquad\text{in distribution}.
\end{align}
\end{theorem}

See \cite{BekermanFigalliGuionnet,BourgadeErdosYau,DeiftGioev,KrishnapurRiderVirag,PasturShcherbina,RamirezRiderVirag,Shcherbina}
for the details, as well as \cite{Bekerman} for a multi-cut extension of Theorem \ref{Theorem: Edge Universality of Beta-Ensembles},
whose statement we omit for simplicity.
Given that $r_0=0$ when $\theta=(1,\be,\infty)$,
we thus obtain the following corollary from Theorem \ref{Theorem: Main Informal}:

\begin{corollary}
\label{Corollary: Beta Ensembles}
Suppose that \eqref{Equation: Beta-Ensembles Scaling Limit} holds,
and let $\mc T$ be the deterministic function in \eqref{Equation: Main Informal}/\eqref{Equation: Main}. For every $\be>0$, one has
\[\mc T\left(\lim_{n\to\infty}\mf s_n\big(\mf e-\mc B^{\be,n}\big)\right)=1/\be
\qquad\text{almost surely},\]
where the limit above is in distribution.
\end{corollary}

Similarly to Theorem \ref{Theorem: Edge Universality of Beta-Ensembles}, there is an obvious
multi-cut extension to Corollary \ref{Corollary: Beta Ensembles}, which we omit for simplicity.

\subsubsection{Open Problem}

Following-up on Corollary \ref{Corollary: Beta Ensembles}, it is natural to wonder
if it is possible to construct consistent estimators of $\be$ that (1) only depend on $\mc B^{\be,n}$'s edge fluctuations
and (2) can be applied to finite-dimensional
models. More specifically:

\begin{openproblem}
\label{Open Problem: Beta Ensembles}
Let $\mf e$ and $\mf s_n$ be as in \eqref{Equation: Beta-Ensembles Scaling Limit}.
Does there exist a sequence of measurable functions $\{\mc T_n:n\in\mbb N\}$ such that,
as $n\to\infty$,
\begin{enumerate}[(1)]
\item $\mc T_n$ only depends on the particles $\mc B^{\be,n}_k\in\mc B^{\be,n}$
such that $|\mc B^{\be,n}_k-\mf e|=o(1)$; and
\item $\mc T_n(\mc B^{\be,n})\to1/\be$ almost surely?
\end{enumerate}
If so, can one also characterize the fluctuations of $\mc T_n(\mc B^{\be,n})-1/\be$ as $n\to\infty$?
\end{openproblem}

Thanks to Corollary \ref{Corollary: Beta Ensembles} and \eqref{Equation: Main}, we expect that a positive
resolution to Open Problem \ref{Open Problem: Beta Ensembles} could potentially rely on establishing a statement of the form
\begin{align}
\label{Equation: Convergence of Finite-Dimensional Traces}
\lim_{n\to\infty}\left|\sum_{k\geq0\text{ s.t. }|\mc B^{\be,n}_k-\mf e|\leq\eps_n}\mr e^{-t_n\,\mf s_n(\mf e-\mc B^{\be,n}_k)/2}-\sum_{k=0}^\infty\mr e^{-t_n\La^{(1,\be,\infty)}_k/2}\right|=0
\end{align}
for $\eps_n,t_n=o(1)$. A similar result already exists in the special case where $V(x)=x^2/2$ when $t_n$ is replaced by a fixed $t>0$;
see the proof of \cite[Corollary 2.10]{GorinShkolnikov}.

\subsection{Detecting Critical Signals in Spiked Models}
\label{Section: Spiked Models}

\subsubsection{Background}

In the early 2000's, Johnstone introduced the spiked Wishart ensemble
\cite[Pages 301--304]{Johnstone} as a means of modeling the effect of
low-rank signals on the spectral statistics of covariance estimators in high dimension.
Some years later, P\'ech\'e \cite{Peche} introduced a closely related additive spiked ensemble to study
low-rank perturbations of self-adjoint matrices.
In their simplest incarnation (i.e., rotationally-symmetric noise),
these models are defined as follows:
Let $r\in\mbb N$, $\be\in\{1,2,4\}$, and $\ell=(\ell_{r-1},\ldots,\ell_0)\in[0,\infty)^r$ be fixed,
assuming that $\ell_{r-1}\leq\cdots\leq\ell_0$.
\begin{enumerate}[(1)]
\item For any integers $n\geq0$ and $p\geq r$, let $D_{\be,n,p}$ be a $n\times p$ matrix whose
entries are i.i.d. standard Gaussians in $\mbb F_\be$, and let $\Si_{\ell,p}$ be the $p\times p$ diagonal matrix with eigenvalues
$1,\ldots,1,1+\ell_{r-1},\ldots,1+\ell_0$. We define the $(n,p)$-dimensional spiked Wishart ensemble as follows ($\cdot^*$ denotes the conjugate transpose):
\[W_{\be,\ell,n,p}=\tfrac1nD_{\be,n,p}\Si_{\ell,p}D_{\be,n,p}^*.\]
\item For any integer $n\geq r,$ let $X_{\be,n}=(D_{\be,n}^*+D_{\be,n})/\sqrt 2$,
where $D_{\be,n}$ is an $n\times n$ matrix whose entries are i.i.d. standard Gaussians in $\mbb F_\be$
(in other words, $X_{\be,n}$ is a GOE, GUE, or GSE matrix, depending on whether $\be=1$, $2$, or $4$).
Then, let $P_{\ell,n}$ be the $n\times n$ diagonal matrix with eigenvalues $0,\ldots,0,\ell_{r-1},\ldots,\ell_0$.
We define the $n$-dimensional spiked Gaussian invariant ensemble as
\[Y_{\be,\ell,n}=\tfrac1{\sqrt n}X_{\be,n}+P_{\ell,n}.\]
\end{enumerate}
The main problem involving these models consists of detecting the signal $\ell$
using the noisy observations provided by the spectra of $W_{\be,\ell,n,p}$ or $Y_{\be,\ell,n}$
in the high-dimensional setting $n\asymp p$.
As a first step toward this goal, one may consider the asymptotic version
of the question:

\begin{question}
\label{Question: Spike Recovery}
Let $r\in\mbb N$ and $\be=\{1,2,4\}$ be fixed.
In the case of the Wishart ensemble, we assume that $p=p(n)$ depends on $n$
in such a way that
\begin{align}
\label{Equation: Wishart Gamma}
p(n)=\ga n\big(1+o(1)\big)\qquad\text{as }n\to\infty
\end{align}
for some $\ga\in(0,\infty)$.
Let $\ell=\ell(n)\in[0,\infty)^r$ be either fixed or dependent on $n$, and let
$\mc S^{\be,\ell,n}=\big\{\mc S^{\be,\ell,n}_1,\ldots,\mc S^{\be,\ell,n}_n\big\}$
be the eigenvalues of either $W_{\be,\ell,n,p(n)}$ or $Y_{\be,\ell,n}$, in increasing order.
Are there asymptotics of $\mc S^{\be,\ell,n}$
(as $n\to\infty$) from which any information about
the signal $\ell$ can be recovered almost surely?
\end{question}

The first papers in this direction were by Baik-Ben Arous-P\'ech\'e
\cite{BaikBenArousPeche} and P\'ech\'e \cite{Peche}. Among other things, these results showed that
the almost-sure asymptotic detectability of each signal $\ell_i$ depends on its size relative to the threshold
\[\tau=\begin{cases}
\sqrt{\ga}&\text{(Wishart)},\\
1&\text{(Gaussian invariant)}.
\end{cases}\]
That is, the analysis of Question \ref{Question: Spike Recovery} is typically split into
three distinct cases:
\begin{enumerate}[(1)]
\item The supercritical case, where $\ell_i>\tau$ for at least one $0\leq i\leq r-1$;
\item the subcritical case, where $\ell_{r-1}\leq\cdots\leq\ell_0<\tau$; and
\item the critical case, where $\ell_{r_{\mr c}-1},\ldots,\ell_0=\tau+o(1)$ for some $r_{\mr c}\geq1$.
\end{enumerate}
To the best of our knowledge, the current state of the art for this theory can be briefly summarized as follows.
Firstly, we state a result for the supercritical case:

\begin{theorem}[\cite{BaikBenArousPeche,BaikSilverstein,CapitaineDonatiMartinFeral,Peche}]
Let $\mf e$ be the top edge of the limiting empirical spectral distribution (ESD) of the
unperturbed Wishart or Gaussian invariant model, i.e.,
\begin{align}
\label{Equation: Spiked Edge}
\mf e=\begin{cases}
(1+\sqrt\ga)^2&\text{(Wishart)},\\
2&\text{(Gaussian invariant)}.
\end{cases}
\end{align}
Define the function
\[\mc O(x)=\begin{cases}
(1+x)(1+\ga/x)&\text{(Wishart)},\\
x+1/x&\text{(Gaussian invariant)}.
\end{cases}\]
For any fixed $\ell\in[0,\infty)^r$ and $i\geq0$, one has
\begin{align}
\label{Equation: Spiked Supercritical Outlier}
\lim_{n\to\infty}\mc S^{\ell,\be,n}_{n-i}=
\begin{cases}
\mc O(\ell_i)&\text{if $i\leq r-1$ and $\ell_i>\tau$},\\
\mf e&\text{if $i\leq r-1$ and $\ell_i\leq\tau$, or if $i>r-1$}
\end{cases}
\qquad\text{almost surely}.
\end{align}
\end{theorem}

\begin{remark}
\label{Remark: Supercritical Only Needs Finite}
$\mc O(x)>\mf e$ is strictly increasing for $x\in(\tau,\infty)$. Thus,
\eqref{Equation: Spiked Supercritical Outlier} implies that the number and magnitude of
all supercritical signals $\ell_i>\tau$ can be recovered from the
asymptotics of a {\it finite} number of the top eigenvalues of $W_{\be,\ell,n,p}$ and $Y_{\be,\ell,n}$.
\end{remark}

Secondly, we state a result regarding the impossibility of detecting subcritical signals
almost-surely using the eigenvalues $\mc S^{\be,\ell,n}$ alone:

\begin{theorem}[\cite{JungChungLee,MontanariReichmanZeitouni,OnatskiMoreiraHallin,OnatskiMoreiraHallin2}]
\label{Theorem: Mutual Contiguity}
Let $\be=1$ (i.e., $\mbb F_\be=\mbb R$).
For any fixed $\ell\in[0,\infty)^r$ such that $\ell_{r-1}\leq\cdots\leq\ell_0<\tau$,
the sequence of spiked eigenvalues $(\mc S^{\be,\ell,n}:n\geq 1)$ is contiguous with respect to the
sequence of unperturbed eigenvalues $(\mc S^{\be,0,n}:n\geq 1)$. That is, for any sequence of
measurable maps $\{\mc T_n:n\geq1\}$,
\[\lim_{n\to\infty}\mbf P\big[\mc T_n(\mc S^{\be,0,n})=0\big]=1
\qquad\text{implies that}\qquad
\lim_{n\to\infty}\mbf P\big[\mc T_n(\mc S^{\be,\ell,n})=0\big]=1.\]
\end{theorem}

\begin{remark}
\label{Remark: Subcritical Uses All}
Despite the impossibility result in Theorem \ref{Theorem: Mutual Contiguity},
one could ask if there nevertheless exist tests that distinguish between the hypotheses
\begin{align}
\label{Equation: Hypotheses}
\bs H_0\text{ : "no perturbation"}
\qquad\text{versus}\qquad
\bs H_\ell\text{ : "perturbation $=\ell$"}
\end{align}
with nontrivial statistical power when $\ell\neq0$ is subcritical.
The short answer to this question is yes; see, e.g., \cite{ElAlaouiKrzakalaJordan,JungChungLee,JohnstoneOnatski,OnatskiMoreiraHallin,OnatskiMoreiraHallin2}
and references therein for a non-exhaustive sample of this very sizeable literature.
A unifying feature of all these results is that the statistical tests that achieve optimal or close-to-optimal power in distinguishing \eqref{Equation: Hypotheses}
(e.g., likelihood ratios or linear spectral statistics)
all use the full spectrum $\mc S^{\be,\ell,n}_1\leq\cdots\leq\mc S^{\be,\ell,n}_n$ as an input. In the words of Dobriban \cite{Dobriban}:
\begin{quote}
"{\it All eigenvalues matter to achieve sharp detection
of weak principal components in high-dimensional data.}"
\end{quote}
To give a simple example of how such results look like, consider the following:
\end{remark}

\begin{theorem}[\cite{ElAlaouiKrzakalaJordan,JungChungLee}]
Suppose that $\be=1$, and that we consider the Gaussian invariant model $Y_{\be,\ell,n}$ (hence $\tau=1$).
Let $\ell\in[0,\infty)^r$ be such that $\ell_{r-1}=\cdots=\ell_0=\la$ for some fixed $\la\in[0,1)$.
Let $\mc L_n$ denote the output of the likelihood ratio test (LRT) for \eqref{Equation: Hypotheses}
in the $n$-dimensional model,
where $\mc L_n(\mc S^{\be,\ell,n})=0$ means accepting $\bs H_0$ and $\mc L_n(\mc S^{\be,\ell,n})=1$
means rejecting $\bs H_0$. Define the error function
\begin{align}
\label{Equation: Finite n Type I and II errors}
\mc E_n(\la)=\mbf P\big[\mc L_n(\mc S^{\be,\ell,n})=1~\big|~\bs H_0\text{ is true}\big]+\mbf P\big[\mc L_n(\mc S^{\be,\ell,n})=0~\big|~\bs H_\ell\text{ is true}\big].
\end{align}
The LRT minimizes \eqref{Equation: Finite n Type I and II errors} over all
measurable functions of $\mc S^{\be,\ell,n}$. Moreover,
\begin{align}
\label{Equation: Limit Gaussian Errors}
\lim_{n\to\infty}\mc E_n(\la)
=\mr{erfc}\left(\tfrac r4\sqrt{-\log(1-\la)}\right)\in(0,1),
\end{align}
where we use $\mr{erfc}(z)=\frac{2}{\sqrt\pi}\int_z^\infty\mr e^{-t^2}\d t$ to denote the complementary error function.
\end{theorem}

Thirdly, we did not find any works in the literature that directly address
Question \ref{Question: Spike Recovery} in the critical case $\ell_{r_{\mr c}-1},\ldots,\ell_0=\tau+o(1)$. However,
as pointed out in \cite[Appendix A]{PerryWeinBandeiraMoitra}, consistent estimators
of $\ell$ in this case can be proved to exist by examining the asymptotic power of subcritical tests when
$\ell_i\to\tau^-$ (e.g., \eqref{Equation: Limit Gaussian Errors} goes to zero when $\ell_i=\la\to1^-$).
Thus, the answer to Question \ref{Question: Spike Recovery} is positive in the critical case.
However, to the best of our knowledge, since this conclusion relies on the optimal subcritical estimators,
it follows from Remark \ref{Remark: Subcritical Uses All} that
{\it all eigenvalues} are used in the currently-known methods to consistently recover critical signals.

\subsubsection{New Result}

Remarks \ref{Remark: Supercritical Only Needs Finite} and
\ref{Remark: Subcritical Uses All} highlight a very sharp dichotomy between
the number of eigenvalues in $\mc S^{\be,\ell,n}$ that are required to optimally infer supercritical and subcritical signals
(respectively, a finite number of top eigenvalues versus
all eigenvalues). Given that critical signals interpolate
between these two extremes, it is thus natural to ask:

\begin{question}
In the context of Question \ref{Question: Spike Recovery},
what is the minimal number of
eigenvalues that are needed in order to detect critical signals with probability one?
\end{question}

According to \eqref{Equation: Spiked Supercritical Outlier}, it is impossible
to detect critical signals using any finite collection of top eigenvalues. In this view, the main
result of this paper implies that critical signals "close enough" (see \eqref{Equation: Close Enough to Critical Threshold})
to the threshold $\tau$ can be recovered from the vanishing proportion of top eigenvalues that
contribute to edge scaling limits.
This is based on the following:

\begin{theorem}
[\cite{BloemendalVirag1,BloemendalVirag2}; see also \cite{BaikBenArousPeche,Mo,Peche,Wang}]
\label{Theorem: Spiked Edge Fluctuations}
Let $r\in\mbb N$ and $\be=\{1,2,4\}$ be fixed,
let $\ga$ and $p=p(n)$ be as in \eqref{Equation: Wishart Gamma} in the Wishart case, and
suppose that $\ell=\ell(n)\in[0,\infty)^r$
is either fixed or depends on $n$.
Suppose that for every $1\leq i\leq r$, the limit
\begin{align}
\label{Equation: wi Assumption}
w_i=\begin{cases}
\displaystyle\tfrac1{\sqrt\ga}\left(1+\tfrac1{\sqrt{\ga}}\right)^{-2/3}\lim_{n\to\infty}n^{1/3}\big(\sqrt\ga-\ell_{i-1}\big)&\text{(Wishart)},\\
\displaystyle\lim_{n\to\infty}n^{1/3}(1-\ell_{i-1})&\text{(Gaussian invariant)},
\end{cases}
\end{align}
exists and is an element of $(-\infty,\infty]$; then let $w=(w_1,\ldots,w_r)$.
If we define
\[\mf s_n=\begin{cases}
\tfrac1{\sqrt\ga}\left(1+\tfrac1{\sqrt{\ga}}\right)^{-4/3}n^{2/3}&\text{(Wishart)},\\
n^{2/3}&\text{(Gaussian invariant)},
\end{cases}\]
and we let $\mf e$ be the top edge of the limiting ESD defined in \eqref{Equation: Spiked Edge}, then it holds that
\begin{align}
\label{Equation: Spiked Scaling Limit}
\lim_{n\to\infty}\mf s_n\big(\mf e-\mc S^{\be,\ell,n}\big)=\La^{(r,\be,w)}\qquad\text{in distribution}.
\end{align}
\end{theorem}

Thanks to \eqref{Equation: wi Assumption}, we see that $w_i\in\mbb R$ corresponds to
\begin{align}
\label{Equation: Close Enough to Critical Threshold}
\ell_{i-1}=\tau+O(n^{-1/3}),
\end{align}
whereas $w_i=\infty$ corresponds to $\ell_{i-1}$ being smaller than $\tau$ by a difference of greater
order than $n^{-1/3}$. Thus, in the present context, the parameter $r_0$ introduced in \eqref{Equation: r0} counts the number
of critical signals that are no farther than $O(n^{-1/3})$ from the threshold $\tau$. Given that
one would presumably know whether $\be=1,2,4$ if presented with a realization of the random matrices $W_{\be,\ell,n,p}$
or $Y_{\be,\ell,n}$, Theorem \ref{Theorem: Main Informal} then implies that the number of critical
signals counted by $r_0$ can be determined from the top edge fluctuations of the random matrices:

\begin{corollary}
\label{Corollary: Spiked}
Let $\mc T$ be as in \eqref{Equation: Main Informal}/\eqref{Equation: Main}.
Under the assumptions of Theorem \ref{Theorem: Spiked Edge Fluctuations},
\[\mc T\left(\lim_{n\to\infty}\mf s_n\big(\mf e-\mc S^{\be,\ell,n}\big)\right)-1/\be=r_0\qquad\text{almost surely},\]
where the limit above is in distribution.
\end{corollary}

\subsubsection{Open Problem}

Similarly to our result on beta-ensembles, the natural follow-up to Corollary \ref{Corollary: Spiked}
is the following:

\begin{openproblem}
\label{Open Problem: Spiked}
Let $\mf e$ and $\mf s_n$ be as in \eqref{Equation: Spiked Scaling Limit}.
Does there exist a sequence of measurable functions $\{\mc T_n:n\in\mbb N\}$ such that,
as $n\to\infty$,
\begin{enumerate}[(1)]
\item $\mc T_n$ only depends on the particles $\mc S^{\be,\ell,n}_k\in\mc S^{\be,\ell,n}$
such that $|\mc S^{\be,\ell,n}_k-\mf e|=o(1)$; and
\item $\mc T_n(\mc S^{\be,\ell,n})\to r_0$ almost surely?
\end{enumerate}
If so, can one also characterize the fluctuations of $\mc T_n(\mc S^{\be,\ell,n})-r_0$ as $n\to\infty$?
\end{openproblem}

That being said, the convergence result analogous to \eqref{Equation: Convergence of Finite-Dimensional Traces} that
one would presumably need to prove Open Problem \ref{Open Problem: Spiked} using
\eqref{Equation: Main Informal}/\eqref{Equation: Main} is much more difficult in the case
of spiked models.
To illustrate this, we note that the limit
\[\lim_{n\to\infty}\sum_{k=1}^n\mr e^{-t\,\mf s_n(\mf e-\mc S^{\be,\ell,n})/2}=\sum_{k=0}^\infty\mr e^{-t\La^{(r,\be,w)}/2}\]
for fixed $t>0$ has not been proved for any spiked model when $r_0>0$, even in the simplest case of the scalar SAO (i.e., $r=r_0=1$).
Indeed, using current techniques, such a result would first require
proving a seemingly difficult strong invariance principle for reflected random walks; see \cite[Conjectures 2.23 and 6.11]{ConvergenceToSpiked}.

\subsection{Rigidity}
\label{Section: Rigidity}

\subsubsection{Background}

In the last 15 years, a sizeable literature on "rigidity" properties of
point processes has been developing, starting with pioneering works by
Ghosh and Peres \cite{Ghosh,GhoshPeres}.
We refer to \cite{Coste,GhoshLebowitz} and references therein for a survey.
In this literature, rigidity refers to the observation that for certain strongly-correlated
point processes, some features of the point configuration inside a bounded set
is almost-surely determined by the configuration of points outside that set.
For instance:

\begin{definition}
\label{Definition: Rigidity}
Let $\mc X\subset\mbb R^d$ be a point process.
We say that $\mc X$ is number rigid if for every bounded Borel set $B\subset\mbb R$,
there exists a measurable map $\mc M_B$ such that
\[|\mc X\cap B|=\mc M_B(\mc X\cap B^c)\qquad\text{almost surely}.\]
(In words, the points outside $B$ determine the number of points inside $B$.)
More generally, given any integer $k\geq0$, we say that $\mc X$ is $k$-moment rigid
if for every bounded Borel set $B\subset\mbb R$,
there exists a measurable map $\mc M_B$ such that
\[\sum_{x\in\mc X\cap B}x^k=\mc M_B(\mc X\cap B^c)\qquad\text{almost surely}.\]
(In words, the points outside $B$ determine the $k^{\mr{th}}$ moment of points inside $B$.)
\end{definition}

Examples of nontrivial rigid point processes include the following:

\begin{example}
\label{Example: Ginibre}
The Ginibre point process was shown to be
number rigid in \cite{GhoshPeres}.
\end{example}

\begin{example}
\label{Example: GAF}
The $\al$-Gaussian analytic function
($\al$-GAF) zero point process
(i.e., the zeros of the
random entire function
$z\mapsto\sum_{k\geq0}\frac{g_k}{(k!)^{\al/2}}z^k$, where $g_k$ are i.i.d.
complex Gaussians) were shown to be $k$-moment rigid for all $0\leq k\leq\lfloor1/\al\rfloor$ in
\cite{GhoshKrishnapur,GhoshPeres}.
\end{example}

\begin{remark}
\label{Remark:Rigidity 1}
Rigidity was in part introduced
to help distinguish structured point processes that are otherwise difficult to
tell apart. As stated in \cite[Page 1794]{GhoshPeres}:
\begin{quote}
"{\it While a simple visual inspection suffices
to (heuristically) distinguish a sample of the Poisson process from that of either the
Ginibre or the GAF zero process (of the same intensity), the latter two are hard to set
apart between themselves. It is therefore an interesting question to devise mathematical statistics that distinguish them.}"
\end{quote}
From this point of view, the notion of rigidity becomes especially interesting
once we note that the results outlined in
Examples \ref{Example: Ginibre} and \ref{Example: GAF}
are optimal, and thus serve as a means of distinguishing
the point processes described therein:
\begin{enumerate}[(1)]
\item For the Ginibre process, it was shown in \cite{GhoshPeres} that
the number of points is the only information about the configuration inside a bounded
set $B$ that can be obtained from the configuration outside $B$.
\item For the $\al$-GAF zeros, it was shown in
\cite{GangopadhyayGhoshTan,GhoshPeres} that the $k^{\mr{th}}$ moments
for $0\leq k\leq\lfloor1/\al\rfloor$ is the only information about the configuration inside a bounded
set $B$ that can be obtained from the configuration outside $B$.
\end{enumerate}
A formal version of these optimality results can be stated using the notion of tolerance; see,
e.g., \cite[Definitions 2 and 3]{GhoshKrishnapur}.
\end{remark}

Going back to our main topic,
the number rigidity of $\La^\theta$ was proved in a sequence of
three papers. Namely, first for $\theta=(1,2,\infty)$ in \cite{Bufetov16}, then for $\theta=(1,\be,w)$ with any choice of $\be>0$ and $w\in(-\infty,\infty]$ in \cite{RigiditySAO},
and finally for all $\theta\in\Theta$ in \cite{RigidityMSAO}.

\subsubsection{New Results}
\label{Section: Rigidity New Results}

Theorem \ref{Theorem: Main Informal} refines current results on $\La^\theta$'s rigidity
in two directions: On the one hand, while the results of \cite{RigidityMSAO,RigiditySAO}
proved the existence of a single function $\mc M$ such that
\begin{align}
\label{Equation: M Function for SAO}
|\La^\theta\cap B|=\mc M(\La^\theta\cap B^c)\qquad\text{almost surely}
\end{align}
for every bounded Borel set $B$,
they did not provide an explicit construction of $\mc M$. In this view,
Theorem \ref{Theorem: Main Informal} improves on the latter
by identifying the exact mechanism that creates number rigidity:
If we restrict the trace in \eqref{Equation: Informal Trace Asymptotics} to the eigenvalues outside of
some bounded set, then the fact that $\mr e^{-t\La^\theta_k/2}=1+o(1)$ almost surely
as $t\to0$ for all fixed $k\geq0$ implies the following:

\begin{corollary}
\label{Corollary: Rigidity Explanation}
For every bounded Borel set $B\subset\mbb R$,
it holds that
\begin{align}
\label{Equation: Rigidity Explanation}
\sum_{k\geq0\text{ s.t. }\La^\theta_k\in B^c}\mr e^{-t\La^\theta_k/2}=\sqrt{\frac2\pi}t^{-3/2}+\frac12\left(r_0+\frac1{\be}\right)-\frac14-|\La^\theta\cap B|+\mf o_t,\qquad t\in(0,1],
\end{align}
where the random remainder terms $\mf o_t$ satisfy \eqref{Equation: Informal Trace Asymptotics Condition 1} and \eqref{Equation: Informal Trace Asymptotics Condition 2}.
This provides a clear blueprint to construct the function $\mc M$ in \eqref{Equation: M Function for SAO} explicitly;
see Remark \ref{Remark: Rigidity Function} for details.
\end{corollary}

On the other hand, while
the statement that $\La^\theta$ is number rigid has intrinsic interest---as it
says something nontrivial about correlations within SAO spectra---it nevertheless falls short
of the ambition for the theory of rigidity laid out in Remark \ref{Remark:Rigidity 1}:
Given that SAO eigenvalue point processes remain very similar as we vary the parameter $\theta$,
we do not expect that different instances of $\La^\theta$ satisfy different notions of rigidity, such as the
rigidity/tolerance hierarchy observed for the Ginibre process and the $\al$-GAF zeroes in
Examples \ref{Example: Ginibre} and \ref{Example: GAF}.

\begin{remark}
To the best of our knowledge, the optimality of number rigidity for $\La^\theta$ (i.e., tolerance)
has not been proved.
That said, we expect that such a result could be proved using the technology developed in \cite{GangopadhyayGhoshTan}.
(The mutual absolute continuity of Palm measures of the Airy process---i.e., $\La^{(1,2,\infty)}$---with the same number
of fixed components proved in \cite{Bufetov18} serves as further evidence of this point.)
\end{remark}

In this context, our main result can be viewed as a first step
toward developing statistics that distinguish between different instances of SAO spectra,
which goes beyond the rigidity/tolerance hierarchy:

\begin{corollary}
\label{Corollary: Parameter Rigidity}
If $\theta=(r,\be,w)$ and $\tilde\theta=(\tilde r,\tilde\be,\tilde w)$ are such that $r_0+1/\be\neq\tilde r_0+1/\tilde{\be}$,
then $\La^\theta$ and $\La^{\tilde\theta}$ are mutually singular.
\end{corollary}

\subsubsection{Open Problems}

Following-up on Corollary \ref{Corollary: Parameter Rigidity},
it is natural to ask if there are other examples of parameter-dependent
point processes that are mutually singular despite being
known/expected to satisfy the same notions of rigidity.
As a first step in this direction, we propose the following:

\begin{openproblem}
Consider the family $\{$Sine$_\be:\be>0\}$, which describes the bulk scaling
limits of various random matrices and beta-ensembles \cite{ValkoVirag,ValkoVirag2}.
These are known to be number rigid and tolerant \cite{ChhaibiNajnudel,DLR}.
Are Sine$_\be$ and Sine$_{\tilde\be}$ mutually singular when $\be\neq\tilde\be$?
If so, can one construct an explicit statistic that allows to recover $\be$ almost surely?
\end{openproblem}

Lastly, one can ask about the optimality of Corollary \ref{Corollary: Parameter Rigidity}:

\begin{openproblem}
Find a function $\Psi$ on the parameter space $\Theta$ such that $\La^\theta$ and $\La^{\tilde\theta}$
are mutually singular if {\it and only if} $\Psi(\theta)\neq\Psi(\tilde\theta)$.
Once $\Psi$ is found, characterize the possible relationships between
the laws of $\La^\theta$
and $\La^{\tilde\theta}$ when $\Psi(\theta)=\Psi(\tilde\theta)$
(e.g., equality, mutual absolute continuity, etc.).
\end{openproblem}

\subsection{Cancellation of Subcritical Coordinates}
\label{Section: Cancellations}

\subsubsection{Background}

In \cite[Page 2734]{BloemendalVirag2}, Bloemendal and Vir\'ag made the curious
observation that any given spiked Wishart or Gaussian invariant model converges
to an infinite number of SAOs with different parameters: Looking back at
Theorem \ref{Theorem: Spiked Edge Fluctuations}, if we add a trivial component
to the spike $\ell=(\ell_{r-1},\ldots,\ell_0)$, thus transforming it into $\tilde\ell=(0,\ell_{r-1},\ldots,\ell_0)$, then the matrix model
remains unchanged. However, the parameter $r$ in the limit \eqref{Equation: Spiked Scaling Limit}
increases by one, and the added component has a Dirichlet boundary condition (since $\ell_{r}=0$ corresponds
to $w_{r+1}=\infty$ in \eqref{Equation: wi Assumption}).
If we combine this observation with the fact that the SAOs are obviously
invariant with respect to permutations of their components (as their noises are induced by GOE,
GUE, and GSE matrices), then we obtain the following from Theorem \ref{Theorem: Spiked Edge Fluctuations}:

\begin{corollary}[\cite{BloemendalVirag2}]
\label{Corollary: Cancellation of Subcritical Coordinates}
Let $\be\in\{1,2,4\}$, and let $\theta=(r,\be,w),\tilde\theta=(\tilde r,\be,\tilde w)\in\Theta$ be such that
\begin{align}
\label{Equation: Cancellation Condition 1}
r_0=\tilde r_0
\end{align}
(i.e., $\mc H_\theta$ and $\mc H_{\tilde\theta}$ have the same number of
Robin boundary conditions) and
\begin{align}
\label{Equation: Cancellation Condition 2}
\{w_1,\ldots,w_r\}\cap\mbb R=\{\tilde w_1,\ldots,\tilde w_{\tilde r}\}\cap\mbb R
\end{align}
(i.e., the Robin boundary conditions match, up to permutation).
Then, $\La^\theta\deq\La^{\tilde\theta}$.
\end{corollary}

Given the somewhat indirect nature of the proof of this result, it is natural to ask:

\begin{question}
\label{Question: Cancellations}
Can one prove Corollary \ref{Corollary: Cancellation of Subcritical Coordinates}
using the intrinsic properties of the SAOs only (i.e., without reference to the matrix models that converge to the SAOs)?
\end{question}

\subsubsection{New Result}

While our results do not provide a complete answer to Question \ref{Question: Cancellations},
they nevertheless shed new light on some of the delicate relationships that must exist between $\mc H_\theta$'s parameters
in order for Corollary \ref{Corollary: Cancellation of Subcritical Coordinates} to hold.
To see this, in what follows we
let $\theta=(r,\be,w)$ with $r\geq1$,
$\be\in\{1,2,4\}$, and $w\in(-\infty,\infty]^r$ be fixed,
and we let $r_0$ be defined as in \eqref{Equation: r0}.
Given $\ka,\si,\upsilon>0$,
define the parameter $\eta=(\ka,\si,\upsilon)$, and define the Schr\"odinger operator
\[\hat H_{\theta,\eta}=-\tfrac12\tfrac{\dd^2}{\dd x^2}+\ka x+\xi,\]
which acts on vector-valued $f\in L^2\big([0,\infty),\mbb F_\be\big)^r$ with boundary conditions
\[\begin{cases}
f(i,0)=0&\text{if }w_i=\infty\text{ (Dirichlet)},\\
f'(i,0)=w_if(i,0)&\text{if }w_i\in\mbb R\text{ (Robin)},
\end{cases}\qquad 1\leq i\leq r,\]
and where $\xi$ is a $r\times r$ matrix white noise whose diagonal components
$\xi_{i,i}$ have variance $\si^2$, and whose off-diagonal components $\xi_{i,j}=\xi_{j,i}^*$ have variance $\upsilon^2$.
The operators $\hat H_{\theta,\eta}$ can be viewed as a generalization of the SAOs, since if we choose the parameter $\eta$
in such a way that $\ka=\ka(r)=\frac r2$, $\si^2=\frac1{\be}$,
and $\upsilon^2=\frac12$, then $2\hat H_{\theta,\eta}=\mc H_{\theta}$.

In this context, our new insight concerning Corollary \ref{Corollary: Cancellation of Subcritical Coordinates}
is that the trace asymptotic for $\mc H_\theta$ in \eqref{Equation: Informal Trace Asymptotics} is a special case of the more general asymptotic
\begin{align}
\label{Equation: Informal Trace Asymptotics 2}
\mr{Tr}\big[\mr e^{-t\hat H_{\theta,\eta}}\big]=\tfrac{r}{\sqrt{2\pi}\ka}t^{-3/2}+\tfrac14\left(2r_0-r+\tfrac{r\si^2}{\ka}+\tfrac{r(r-1)\upsilon^2}{\ka}\right)+\mf o_t\qquad\text{for all }t\in(0,1],
\end{align}
where $\mf o_t$ satisfy \eqref{Equation: Informal Trace Asymptotics Condition 1} and \eqref{Equation: Informal Trace Asymptotics Condition 2}.
With this in hand, we can obtain necessary conditions on the parameter $\eta$ in order for the
cancellations in Corollary \ref{Corollary: Cancellation of Subcritical Coordinates} to occur.

More specifically, let $\tilde\theta=(\tilde r,\be,\tilde w)$ be related to $\theta$ in such a way that
\eqref{Equation: Cancellation Condition 1} and \eqref{Equation: Cancellation Condition 2} hold.
Then, consider $\eta=(\ka,\si,\upsilon)$ such that, similarly to the SAOs,
\begin{align}
\label{Equation: eta 1}
\ka=\ka(r)\text{ possibly depends on }r,
\end{align}
and conversely for the parameters that determine the noise's variance,
\begin{align}
\label{Equation: eta 2}
\si\text{ and }\upsilon\text{ are independent of }r.
\end{align}
In order for the operators $\hat H_{\theta,\eta}$ and $\hat H_{\tilde\theta,\eta}$ to generate the same eigenvalue
point process,
it is necessary that every term in the asymptotic \eqref{Equation: Informal Trace Asymptotics 2} be independent of $r$
(as this could be different from $\tilde r$),
and instead only depend on the common parameter $r_0=\tilde r_0$.
In order for the leading order term in \eqref{Equation: Informal Trace Asymptotics 2} to be independent of $r$, it is necessary that
$\ka=\ka(r)=\ka_0r$ for some $\ka_0>0$ independent of $r$.
If such is the case, then \eqref{Equation: Informal Trace Asymptotics 2} becomes
\begin{align}
\label{Equation: Informal Trace Asymptotics 2.1}
\mr{Tr}\big[\mr e^{-t\hat H_{\theta,\eta}}\big]=\tfrac{1}{\sqrt{2\pi}\ka_0}t^{-3/2}+\tfrac14\left(2r_0-r+\tfrac{\si^2}{\ka_0}+\tfrac{(r-1)\upsilon^2}{\ka_0}\right)+\mf o_t.
\end{align}
Once this is applied, the only $r$-dependence in \eqref{Equation: Informal Trace Asymptotics 2.1} lies in the second order term,
via the function $r\mapsto(\upsilon^2/\ka_0-1)r$. Since $\upsilon$ is assumed independent of $r$, the only way for this to not depend on $r$ is to have $\upsilon^2=\ka_0$,
in which case \eqref{Equation: Informal Trace Asymptotics 2.1} becomes
\begin{align}
\label{Equation: Informal Trace Asymptotics 2.2}
\mr{Tr}\big[\mr e^{-t\hat H_{\theta,\eta}}\big]=\tfrac{1}{\sqrt{2\pi}\ka_0}t^{-3/2}+\tfrac14\left(2r_0-1+\tfrac{\si^2}{\ka_0}\right)+\mf o_t.
\end{align}
We thus obtain the following result:

\begin{corollary}
\label{Corollary: Cancellations}
Suppose that $\theta$ and $\tilde\theta$ satisfy \eqref{Equation: Cancellation Condition 1} and \eqref{Equation: Cancellation Condition 2},
and that $\eta$ satisfies \eqref{Equation: eta 1} and \eqref{Equation: eta 2}.
If $\hat H_{\theta,\eta}$ and $\hat H_{\tilde\theta,\eta}$ generate the same eigenvalue
point process, then $\ka(r)=r\upsilon^2$. In particular, if $\upsilon^2=\frac12$ (as is the case for the SAOs),
then the cancellation in Corollary \ref{Corollary: Cancellation of Subcritical Coordinates} can only occur if $\ka(r)=r/2$.
\end{corollary}

\subsubsection{Open Problem}

Based on this partial result, we expect that a more detailed study of the joint moments of $\mr{Tr}\big[\mr e^{-t\mc H_\theta/2}\big]$
(over multiple $t$'s) could be used to completely answer Question \ref{Question: Cancellations}.
As this lies outside the scope of the main applications of this paper, we leave this question open
for future investigations.

\section{Trace Asymptotics}
\label{Section: Proof}

We now embark on the proof of our trace
asymptotic. The most general incarnation
of the asymptotic is \eqref{Equation: Informal Trace Asymptotics 2}.
For this purpose, we
recall that,
given $\theta=(r,\be,w)\in\Theta$ and $\eta=(\ka,\si,\upsilon)\in(0,\infty)^3$,
we denote the operator
\[\hat H_{\theta,\eta}=-\tfrac12\tfrac{\dd^2}{\dd x^2}+\ka x+\xi,\]
which acts on vector-valued $f\in L^2\big([0,\infty),\mbb F_\be\big)^r$ with boundary conditions
\[\begin{cases}
f(i,0)=0&\text{if }w_i=\infty\text{ (Dirichlet)},\\
f'(i,0)=w_if(i,0)&\text{if }w_i\in\mbb R\text{ (Robin)},
\end{cases}\qquad 1\leq i\leq r,\]
and where $\xi$ is a $r\times r$ matrix white noise whose diagonal components
$\xi_{i,i}$ have variance $\si^2$, and whose off-diagonal components $\xi_{i,j}=\xi_{j,i}^*$ have variance $\upsilon^2$.
The statement in \eqref{Equation: Informal Trace Asymptotics 2}
can be split into two claims:

\begin{theorem}
\label{Theorem: Main Technical}
For every
$\theta\in\Theta$ and $\eta\in(0,\infty)^3$, one has
\begin{align}
\label{Equation: Main Technical}
\mbf E\left[\mr{Tr}\big[\mr e^{-t\hat H_{\theta,\eta}}\big]\right]=\tfrac{r}{\sqrt{2\pi}\ka}t^{-3/2}+\tfrac14\left(2r_0-r+\tfrac{r\si^2}{\ka}+\tfrac{r(r-1)\upsilon^2}{\ka}\right)+o(1)\qquad\text{as }t\to0^+,
\end{align}
and there exists a constant $C>0$ such that
\begin{align}
\label{Equation: Main Technical 2}
\mbf{Cov}\left[\mr{Tr}\big[\mr e^{-s\hat H_{\theta,\eta}}\big],\mr{Tr}\big[\mr e^{-t\hat H_{\theta,\eta}}\big]\right]\leq C\left(\frac{\min\{s,t\}}{\max\{s,t\}}\right)^{1/4}
\qquad\text{ for every }s,t\in(0,1].
\end{align}
\end{theorem}

The covariance bound in Theorem \ref{Theorem: Main Technical} was proved in \cite[(8.4)]{RigidityMSAO}.
Consequently, we only need to prove \eqref{Equation: Main Technical}.
The remainder of this section is organized as follows:
Before we prove \eqref{Equation: Main Technical},
in Section \ref{Section: Parameter Recovery} we show how Theorem \ref{Theorem: Main Technical} is used to
construct the function $\mc T$ in \eqref{Equation: Main Informal}, and thus prove Theorem \ref{Theorem: Main Informal}.
Then, in Section \ref{Section: Trace Formula},
we state a Feynman-Kac formula proved in \cite{RigidityMSAO} for the expected trace $\mbf E\left[\mr{Tr}\big[\mr e^{-t\hat H_{\theta,\eta}}\big]\right]$.
Finally, in Sections \ref{Section: Main Proof 1} and \ref{Section: Bridge Asymptotics},
we use the Feynman-Kac formula in question to prove the asymptotics in \eqref{Equation: Main Technical}.

\subsection{Proof of Theorem \ref{Theorem: Main Informal}}
\label{Section: Parameter Recovery}

We begin by defining $\mc T$:

\begin{definition}
\label{Definition: Main}
Let $c_1,c_2>0$ be any fixed constants.
Define the vanishing sequence of positive numbers $1\geq t_1>t_2>\cdots>0$ as
\begin{align}
\label{Equation: Main tn}
t_n=\mr e^{-c_1n},\qquad n\geq1,
\end{align}
and let $1\leq N_1<N_2<\cdots$ be the sequence of integers
\begin{align}
\label{Equation: Main Nm}
N_m=\lceil m^{1+c_2}\rceil,\qquad m\geq1.
\end{align}
Given any point configuration $\mc X=\{x_0,x_1,x_2,\ldots\}$ on $\mbb R$, we define
\begin{align}
\label{Equation: Main}
\mc T(\mc X)=\begin{cases}
\displaystyle\frac12+\lim_{m\to\infty}\frac2{N_m}\sum_{n=1}^{N_m}\left(\sum_{k=0}^\infty\mr e^{-t_nx_k/2}-\sqrt{\tfrac2{\pi}}t_n^{-3/2}\right)&\text{if the limit exists,}\\
0&\text{otherwise.}
\end{cases}
\end{align}
\end{definition}

\begin{remark}
There is some degree of arbitrariness in how $t_n$ and $N_m$ were chosen in \eqref{Equation: Main tn}
and \eqref{Equation: Main Nm}. As will be made clear by the proof in this section, any sequence of $t_n$'s
and $N_m$'s that respectively vanish and diverge "fast enough" suffice to obtain \eqref{Equation: Main Informal}.
Our specific choices in \eqref{Equation: Main tn} and \eqref{Equation: Main Nm} are made for simplicity
and to emphasize that it is possible to be completely explicit in our definition of the map $\mc T$.
\end{remark}

Recall that if $\eta=(\frac r2,\frac1\be,\frac12)$, then $2\hat H_{\theta,\eta}=\mc H_{\theta}$. Thus,
\eqref{Equation: Main Technical} implies that
\[\mbf E\left[\mr{Tr}\big[\mr e^{-t\mc H_{\theta}/2}\big]\right]=\sqrt{\tfrac{2}{\pi}}t^{-3/2}+\tfrac{r_0+1/\be}{2}-\tfrac14+o(1)\qquad\text{as }t\to0^+.\]
By a trivial rearrangement of the terms in \eqref{Equation: Main},
we therefore obtain that $\mc T(\La^\theta)=r_0+1/\be$ almost surely if we can prove that
\[\lim_{m\to\infty}\frac{1}{N_m}\sum_{n=1}^{N_m}\left(\mr{Tr}\big[\mr e^{-t_n\mc H_\theta/2}\big]-\mbf E\left[\mr{Tr}\big[\mr e^{-t_n\mc H_{\theta}/2}\big]\right]\right)=0
\qquad\text{almost surely.}\]
Toward this end, by Markov's inequality, we get that for every $m\in\mbb N$ and $\eps>0$,
\begin{align}
\nonumber
&\mbf P\left[\left|\frac{1}{N_m}\sum_{n=1}^{N_m}\left(\mr{Tr}\big[\mr e^{-t_n\mc H_\theta/2}\big]-\mbf E\left[\mr{Tr}\big[\mr e^{-t_n\mc H_{\theta}/2}\big]\right]\right)\right|\geq m^{-\eps/2}\right]\\
\nonumber
\leq&m^\eps~\mbf E\left[\frac1{N_m^2}\left(\sum_{n=1}^{N_m}\mr{Tr}\big[\mr e^{-t_n\mc H_\theta/2}\big]-\mbf E\left[\mr{Tr}\big[\mr e^{-t_n\mc H_\theta/2}\big]\right]\right)^2\right]\\
\label{Equation: Main 2}
=&\frac{m^\epsilon}{N_m^2}\sum_{n_1,n_2=1}^{N_m}\mbf{Cov}\left[\mr{Tr}\big[\mr e^{-t_{n_1}\mc H_\theta/2}\big],\mr{Tr}\big[\mr e^{-t_{n_2}\mc H_\theta/2}\big]\right].
\end{align}
By the Borel-Cantelli lemma, it suffices to show that \eqref{Equation: Main 2} is summable in $m$.

By \eqref{Equation: Main Technical 2}, there is a constant $C>0$ independent of $N$ such that
\[\frac{m^\epsilon}{N_m^2}\sum_{n_1,n_2=1}^{N_m}\mbf{Cov}\left[\mr{Tr}\big[\mr e^{-t_{n_1}\mc H_\theta/2}\big],\mr{Tr}\big[\mr e^{-t_{n_2}\mc H_\theta/2}\big]\right]\leq\frac{Cm^\eps}{N_m^2}\sum_{n_1,n_2=1}^{N_m}\left(\frac{\min\{t_{n_1},t_{n_2}\}}{\max\{t_{n_1},t_{n_2}\}}\right)^{1/4}.\]
With $t_n$ defined as in \eqref{Equation: Main tn}, this simplifies to
\[\frac{Cm^\eps}{N_m^2}\sum_{n_1,n_2=1}^{N_m}\left(\frac{\min\{t_{n_1},t_{n_2}\}}{\max\{t_{n_1},t_{n_2}\}}\right)^{1/4}=\frac{Cm^\eps}{N_m^2}\sum_{n_1,n_2=1}^{N_m}\mr e^{-c_1|n_1-n_2|/4}.\]
We can split the above sum as follows:
\[\frac{Cm^\eps}{N_m^2}\sum_{\substack{1\leq n_1,n_2\leq N_m\\|n_1-n_2|\leq(4/c_1)\log N_m}}\mr e^{-c_1|n_1-n_2|/4}+
\frac{Cm^\eps}{N_m^2}\sum_{\substack{1\leq n_1,n_2\leq N_m\\|n_1-n_2|>(4/c_1)\log N_m}}\mr e^{-c_1|n_1-n_2|/4}.\]
Given that $m^{1+c_2}\leq N_m\leq m^{1+c_2}+1$ (as per \eqref{Equation: Main Nm}) and $\mr e^{-c_1|n_1-n_2|/4}\leq1$, we get
\[\frac{Cm^\eps}{N_m^2}\sum_{\substack{1\leq n_1,n_2\leq N_m\\|n_1-n_2|\leq(4/c_1)\log N_m}}\mr e^{-c_1|n_1-n_2|/4}=O\left(\frac{m^\eps}{N_m^2}\cdot N_m\log N_m\right)=O(m^{-1-(c_2-\eps)}\log m)\]
as $m\to\infty$.
Moreover,
\begin{multline*}
\frac{Cm^\eps}{N_m^2}\sum_{\substack{1\leq n_1,n_2\leq N_m\\|n_1-n_2|>(4/c_1)\log N_m}}\mr e^{-c_1|n_1-n_2|/4}\\
\leq
\frac{Cm^\eps\mr e^{-\log N_m}}{N_m^2}\sum_{n_1,n_2=1}^{N_m}1=O(m^\eps N_m^{-1})=O(m^{-1-(c_2-\eps)})
\end{multline*}
as $m\to\infty$.
If we choose $\eps<c_2$, then this implies that \eqref{Equation: Main 2} is summable,
which concludes the proof that $\mc T(\La^\theta)=r_0+1/\be$ almost surely.

\begin{remark}
\label{Remark: Rigidity Function}
Following-up on the calculations performed in this section, we see that
\[\lim_{m\to\infty}\frac1{N_m}\sum_{n=1}^{N_m}\left(\sqrt{\tfrac2\pi}t_n^{-3/2}+\tfrac12\left(r_0+\tfrac1{\be}\right)-\tfrac14-\sum_{k=0}^\infty\mr e^{-t_n\La^\theta_k/2}\right)=0\qquad\text{almost surely},\]
where $t_n$ and $N_m$ are defined as in \eqref{Equation: Main tn} and \eqref{Equation: Main Nm} respectively.
Moreover, for any $k\geq0$, one has
\[\lim_{m\to\infty}\frac1{N_m}\sum_{n=1}^{N_m}\mr e^{-t_n\La^\theta_k/2}=1\qquad\text{almost surely.}\]
Consequently,
if we define the map
\[\mc M(\mc X)=\begin{cases}
\displaystyle\lim_{m\to\infty}\frac1{N_m}\sum_{n=1}^{N_m}\left(\sqrt{\tfrac2\pi}t_n^{-3/2}+\tfrac12\left(r_0+\tfrac1{\be}\right)-\tfrac14-\sum_{k=0}^\infty\mr e^{-t_nx_k/2}\right)&\text{if the limit exists,}\\
0&\text{otherwise,}
\end{cases}\]
for any point configuration $\mc X=\{x_0,x_1,x_2,\ldots\}\subset\mbb R$,
then \eqref{Equation: M Function for SAO} holds. Indeed, for any bounded set $B\subset\mbb R$,
we can write
\begin{multline}
\label{Equation: Rigidity Proof}
\mc M(\La^\theta\cap B^c)=\lim_{m\to\infty}\frac1{N_m}\sum_{n=1}^{N_m}\left(\sqrt{\tfrac2\pi}t_n^{-3/2}+\tfrac12\left(r_0+\tfrac1{\be}\right)-\tfrac14-\sum_{k:\La^\theta_k\in B^c}^\infty\mr e^{-t_n\La^\theta_k/2}\right)\\
=\lim_{m\to\infty}\frac1{N_m}\sum_{n=1}^{N_m}\left(\sqrt{\tfrac2\pi}t_n^{-3/2}+\tfrac12\left(r_0+\tfrac1{\be}\right)-\tfrac14-\sum_{k=0}^\infty\mr e^{-t_n\La^\theta_k/2}\right)+
\lim_{m\to\infty}\frac1{N_m}\sum_{n=1}^{N_m}\sum_{k:\La^\theta_k\in B}\mr e^{-t_n\La^\theta_k/2};
\end{multline}
the first limit on the second line of \eqref{Equation: Rigidity Proof} is almost-surely zero,
and the second is almost-surely $|\La^\theta\cap B|$.
\end{remark}

\subsection{Proof of \eqref{Equation: Main Technical} Step 1. Expected Trace Formula}
\label{Section: Trace Formula}

The first step in our proof of \eqref{Equation: Main Technical} is to
state the Feynman-Kac formula for $\mbf E\left[\mr{Tr}\big[\mr e^{-t\hat H_{\theta,\eta}}\big]\right]$
obtained in \cite[Theorem 4.9]{RigidityMSAO}, which is what we will use to calculate
the leading orders of the expected trace.
This is stated as Lemma \ref{Lemma: Raw Feynman-Kac} below.
Due to the complexity of the formula in question,
we must set up a sizeable amount of notation to state it. First,
we introduce the main stochastic process that drives
the Feynman-Kac formula, which is the reflected Brownian bridge and its local times:

\begin{definition}
\label{Definition: B and X}
Let $B$ be a standard Brownian motion on $\mbb R$,
and denote $B$'s transition kernel by
\[\Pi_B(t;x,y)=\frac{\mr e^{-(x-y)^2/2t}}{\sqrt{2\pi t}},\qquad t>0,~x,y\in\mbb R.\]
We denote
$B^x=(B|B(0)=x)$ and $B^{x,y}_t=(B|B(0)=x\text{ and }B(t)=y)$.
Then, let $X=|B|$ be a reflected standard Brownian motion on $[0,\infty)$,
and denote $X$'s transition kernel by
\[\Pi_X(t;x,y)=\frac{\mr e^{-(x-y)^2/2t}}{\sqrt{2\pi t}}+\frac{\mr e^{-(x+y)^2/2t}}{\sqrt{2\pi t}},\qquad t>0,~x,y\geq0.\]
We similarly denote
$X^x=(X|X(0)=x)$ and $X^{x,y}_t=(X|X(0)=x\text{ and }X(t)=y)$.
\end{definition}

\begin{definition}
Fix $t>0$,
and let $Z$ be either $B$, $X$, or one of its conditioned version (i.e., $B^x$, $X^x$, $B^{x,y}_t$, or $X^{x,y}_t$).
Given $0\leq u<v\leq t$,
we let $\mf L_{[u,v)}^0(Z)$ denote $Z$'s boundary local time collected on
the time interval $[u,v)$, that is,
\begin{align*}
\mf L^0_{[u,v)}(Z)=\lim_{\eps\to0}\frac1{2\eps}\int_u^v\mbf 1_{\{-\eps<Z(s)<\eps\}}\d s\qquad\text{almost surely}.
\end{align*}
Moreover, we let $\big(L^y_{[u,v)}(Z):y\in[0,\infty)\big)$ denote the continuous version of
$Z$'s local time process on $[u,v)$, that is the process such that
\[\int_u^vf\big(Z(s)\big)\d s=\int_0^\infty L^y_{[u,v)}(Z)f(y)\d y\qquad\text{for every measurable }f:\mbb R\to\mbb R.\]
Finally, we use the shorthands $\mf L_t^0(Z)=\mf L_{[0,t)}^0(Z)$
and $L_t^y(Z)=L_{[0,t)}^y(Z)$.
\end{definition}

If we were only dealing with scalar-valued Schr\"odinger operators on $[0,\infty)$,
then these would be the only processes that are needed to state the Feynman-Kac
formula. However, since we are dealing with the vector-valued setting, we need
more tools. First, we set up some notation to describe the function spaces related
to $\hat H_{\theta,\eta}$:

\begin{definition}
Given $f,g:[0,\infty)\to\mbb F_\be$, we denote
\[\langle f,g\rangle=\int_0^\infty f(x)g(x)\d x\qquad\text{and}\qquad\|f\|_2=\sqrt{\langle f,f\rangle}.\]
We denote $\mc A=\{1,\ldots,r\}\times [0,\infty)$, and we let $\mu$ denote
the measure on $\mc A$ obtained by taking the product of the counting measure
on $\{1,\ldots, r\}$ and the Lebesgue measure on $[0,\infty)$.
Given $f:\mc A\to\mbb F_\be$, we denote
\[\|f\|_\mu=\sqrt{\int_\mc A|f(a)|^2\d\mu(a)}=\sqrt{\sum_{i=1}^r\|f(i,\cdot)\|_2^2}.\]
\end{definition}

Next, the vector-valued setting requires the introduction of a jump process
on $\{1,\ldots,r\}$ which depends on $X$'s self-intersections as follows:

\begin{definition}
\label{Definition: Pn}
Let $\mc P_0=\varnothing$ be the empty set, and
given an even integer $n\geq2$, let $\mc P_n$ denote the set of perfect pair matchings
of $\{1,\ldots,n\}$. That is, $\mc P_n$ is the set
\[\Big\{\big\{\{\ell_1,\ell_2\},\{\ell_3,\ell_4\},\ldots,\{\ell_{n-1},\ell_n\}\big\}:
\ell_i\in\{1,\ldots,n\}\text{ and }
\{\ell_i,\ell_{i+1}\}\cap\{\ell_j,\ell_{j+1}\}=\varnothing\text{ if }j\neq i\Big\}.\]
\end{definition}

\begin{definition}
Let $t>0$ and $x\geq0$, let $n\geq2$ be an even integer, and let $q\in\mc P_n$.
Given a realization of $X^{x,x}_t$, we let $\mr{si}_{t,n,q,X^{x,x}_t}$ denote
the unique random Borel probability measure on $[0,t)^n$ such that for every
$s_1,t_1,\ldots,s_n,t_n\in[0,t)$ with $s_k<t_k$ for $k=1,\ldots,n$, one has
\begin{multline}
\label{Equation: SI Measures}
\mr{si}_{t,n,q,X^{x,x}_t}\big([s_1,t_1)\times[s_2,t_2)\times\cdots\times[s_n,t_n)\big)\\
=\frac{1}{\|L_t(X^{x,x}_t)\|_2^{2n}}\prod_{\{\ell_1,\ell_2\}\in q}\left\langle L_{[s_{\ell_1},t_{\ell_1})}(X^{x,x}_t),L_{[s_{\ell_2},t_{\ell_2})}(X^{x,x}_t)\right\rangle
\end{multline}
(this uniquely determines the measure since product rectangles are a determining class
for the Borel $\si$-algebra).
\end{definition}

\begin{definition}
\label{Definition: Singular Process}
Let $t>0$ and $x\geq0$ be fixed,
and suppose that we are given a realization of $X^{x,x}_t$.
For any $1\leq i\leq r$, we generate the jump process $\hat U^i_t:[0,t]\to\{1,\ldots,r\}$
using the following step-by-step procedure:
\begin{enumerate}[$\bullet$]
\item Given $X^{x,x}_t$, sample $\hat N(t)$ according to the distribution
\begin{align}
\label{Equation: Hat N Probability}
\mbf P\big[\hat N(t)=2m\big]=\left(\frac{(r-1)^2\|L_{t}(X^{x,x}_t)\|^2_2}{2}\right)^m\frac{1}{m!}\mr e^{-\frac{(r-1)^2}2\|L_{t}(X^{x,x}_t)\|_2^2},\qquad m=0,1,2,\ldots,
\end{align}
i.e., $\hat N(t)/2$ is Poisson with parameter $\frac{(r-1)^2}2\|L_{t}(X^{x,x}_t)\|_2^2$.
\item Given $\hat N(t)$, sample $\hat q$ uniformly at random in $\mc P_{\hat N(t)}$.
\item Given $X^{x,x}_t$, $\hat N(t)$, and $\hat q$,
draw $\bs T=(T_1,\ldots,T_{\hat N(t)})\in[0,t]^{\hat N(t)}$
according to the self-intersection measure $\mr{si}_{t,\hat N(t),\hat q,X^{x,x}_t}$. Then,
we let $\hat{\bs\tau}=(\hat\tau_1,\ldots,\hat\tau_{\hat N(t)})$ be the tuple $\bs T$ rearranged in increasing order.
\item Given $\hat N(t)$, $\bs T$, and $\hat{\bs\tau}$, we let $\hat\pi$ be the random permutation
of $\{1,\ldots,\hat N(t)\}$
such that $\hat\tau_{\hat\pi(\ell)}=T_\ell$ for all $\ell\leq\hat N(t)$, and we let $\hat p\in\mc P_{\hat N(t)}$
be the random perfect matching such that $\{\ell_1,\ell_2\}\in\hat q$ if and only if
$\{\hat\pi(\ell_1),\hat\pi(\ell_2)\}\in\hat p$. In words, $\hat p$ ensures that each pair of times $T_k$
that were matched by $\hat q$ are still matched once arranged in nondecreasing order $\hat \tau_k$.
\item Let $M^i=(M^i_0,M^i_1,M^i_2,\ldots)$ be a uniform random walk on
the complete graph (without self-edges) on $\{1,\ldots,r\}$ with starting point $M^i_0=i$, i.e.,
\begin{align}
\label{Equation: M Probability}
\mbf P\big[M^i_k=i_2~\big|~M^i_{k-1}=i_1\big]=\frac{\mbf 1_{\{i_2\neq i_1\}}}{r-1},\qquad k\geq1,~i_1,i_2\in\{1,\ldots,r\}.
\end{align}
We assume $M^i$
is independent of $X^{x,x}_t$, $\hat N(t)$, and $\hat{\bs\tau}$.
\item Combining all of the above, we then let
$\big(\hat U^i_t(s):1\leq s\leq t\big)$ be the c\`{a}dl\`{a}g
path with jump times $\hat{\bs\tau}$
and jumps $M^i$:
\begin{align}
\label{Equation: U's path in terms of M}
\hat U^i_t(s)=i\mbf 1_{[0,\hat\tau_1)}(s)+\sum_{k=1}^{\hat N(t)-1} M^i_k\mbf 1_{[\hat\tau_k,\hat\tau_{k+1})}(s)+M_{\hat N(t)}\mbf 1_{[\hat\tau_{\hat N(t)},t]}(s),\qquad 0\leq s\leq t.
\end{align}
In particular, if $\hat N(t)=0$, then $\hat U^i_t(t)=i$ is constant.
\end{enumerate}
Then, we define local times of the combined processes $(\hat U^i_t,X^{x,x}_t)$ as follows:
For any $t>0$ and $(j,y)\in\mc A=\{1,\ldots,r\}\times [0,\infty)$, we let
\begin{align}
\label{Equation: Combined Local Time 1}
\mf L^{(j,0)}_t(\hat U^i_t,X^{x,x}_t)
=\sum_{k\geq1:~U^i_t(\hat\tau_{k-1})=j}\mf L_{[\hat\tau_{k-1},\hat\tau_k)\cap[0,t]}^0(X^{x,x}_t)
\end{align}
and
\begin{align}
\label{Equation: Combined Local Time 2}
L^{(j,y)}_t(\hat U^i_t,X^{x,x}_t)
=\sum_{k\geq1:~U^i_t(\hat\tau_{k-1})=j}L_{[\hat\tau_{k-1},\hat\tau_k)\cap[0,t]}^y(X^{x,x}_t)
\end{align}
respectively denote $X^{x,x}_t$'s boundary local time at zero and local time
collected on the intervals of $[0,t]$ where the jump process $\hat U^i_t$ is equal to $j$.
We note that since $[u,v)\mapsto\mf L^0_{[u,v)}(X^{x,x}_t)$ and $[u,v)\mapsto L^y_{[u,v)}(X^{x,x}_t)$
induce measures on $[0,\infty)$, it is clear that
\begin{align}
\label{Equation: Combined Local Time Sum}
\sum_{j=1}^r\mf L^{(j,0)}_t(\hat U^i_t,X^{x,x}_t)=\mf L_t(X^{x,x}_t)
\quad\text{and}\quad
\sum_{j=1}^rL^{(j,y)}_t(\hat U^i_t,X^{x,x}_t)=L^y_t(X^{x,x}_t).
\end{align}
\end{definition}

Finally, the computation of the moment $\mbf E\Big[\mr{Tr}\big[\mr e^{-t\hat H_{\theta,\eta}}\big]\Big]$
requires introducing the following combinatorial constant:

\begin{definition}
\label{Definition: Bn}
Given $n\in\mbb N$, let
$\mc B_n$ denote the set of binary sequences with $n$ steps,
namely, sequences of the form $m=(m_0,m_1,\ldots,m_n)\in\{0,1\}^{n+1}$.
Given $h,l\in\{0,1\}$, we let $\mc B^{h,l}_n$ denote the set of
sequences $m\in\mc B_n$ such that $m_0=h$ and $m_n=l$.
\end{definition}

\begin{definition}
\label{Definition: Combinatorial Constant}
Let $t>0$ and $1\leq i\leq r$ be fixed.
Suppose we are given a realization of $\hat N(t)\geq2$, the matching $\hat p\in\mc P_{\hat N(t)}$,
and the process $\hat U^i_t$
in Definition \ref{Definition: Singular Process}.
Letting $M^i=(M^i_0,M^i_1,\ldots)$ be as in Definition \ref{Definition: Singular Process}
(i.e., the process that determines the steps taken by $\hat U^i_t$
in $\{1,\ldots,r\}$), we denote the jumps
\[J_k=\big(M^i_{k-1},M^i_k\big),\qquad k\geq1.\]
Then, we let $J_k^*=\big(M^i_k,M^i_{k-1}\big)$ denote these same jumps in reverse order.
We define the combinatorial constant $\mf C_t(\hat p,\hat U^i_t)$ as follows:
\begin{enumerate}[\qquad(a)]
\item If $\mbb F_\be=\mbb R$ (i.e., $\be=1$), then
\[\mf C_t(\hat p,\hat U^i_t)
=\begin{cases}
1&\text{if either $J_{\ell_1}=J_{\ell_2}$ or $J_{\ell_1}=J_{\ell_2}^*$ for all $\{\ell_1,\ell_2\}\in\hat p$},\\
0&\text{otherwise}.
\end{cases}\]
\item If $\mbb F_\be=\mbb C$ (i.e., $\be=2$), then
\[\mf C_t(\hat p,\hat U^i_t)
=\begin{cases}
1&\text{if $J_{\ell_1}=J_{\ell_2}^*$ for all $\{\ell_1,\ell_2\}\in\hat p$},\\
0&\text{otherwise}.
\end{cases}\]
\item If $\mbb F_\be=\mbb H$ (i.e., $\be=4$), then
\[\mf C_t(\hat p,\hat U^i_t)
=\begin{cases}
\mf D_t(\hat p,\hat U^i_t)&\text{if either $J_{\ell_1}=J_{\ell_2}$ or $J_{\ell_1}=J_{\ell_2}^*$ for all $\{\ell_1,\ell_2\}\in\hat p$},\\
0&\text{otherwise}.
\end{cases}\]
Here, we define $\mf D_t(\hat p,\hat U^i_t)$ as follows: A binary sequence $m\in\mc B^{0,0}_{N(t)}$
is said to respect $(\hat p,J)$ if for every
pair $\{\ell_1,\ell_2\}\in\hat p$, the following holds:
\begin{enumerate}[(c.1)]
\item if $J_{\ell_1}=J_{\ell_2}$, then $\big((m_{\ell_1-1},m_{\ell_1}),(m_{\ell_2-1},m_{\ell_2})\big)$
are equal to one of the following pairs: $\big((0,0),(1,1)\big)$, $\big((1,1),(0,0)\big)$,
$\big((0,1),(1,0)\big)$, or $\big((1,0),(0,1)\big)$.
\item if $J_{\ell_1}=J_{\ell_2}^*$, then $\big((m_{\ell_1-1},m_{\ell_1}),(m_{\ell_2-1},m_{\ell_2})\big)$
are equal to one of the following pairs: $\big((0,0),(0,0)\big)$, $\big((1,1),(1,1)\big)$,
$\big((0,1),(1,0)\big)$, or $\big((1,0),(0,1)\big)$.
\end{enumerate}
Lastly, we call any pair $\{\ell_1,\ell_2\}\in p$ such that $J_{\ell_1}=J_{\ell_2}$
and
\[\big((m_{\ell_1-1},m_{\ell_1}),(m_{\ell_2-1},m_{\ell_2})\big)=\big((0,1),(1,0)\big)\text{ or }\big((1,0),(0,1)\big)\]
a flip, and we let $f(m,\hat p,J)$ denote the total number of flips that occur
in the combination of $m,\hat p,$ and $J$.
With all this in hand, we finally define
\begin{align}
\label{Equation: D constant for H}
\mf D_t(\hat p,\hat U^i_t)=2^{-\hat N(t)/2}\sum_{m\in\mc B^{0,0}_{N(t)},~m\text{ respects }(\hat p,J)}(-1)^{f(m,\hat p,J)}.
\end{align}
\end{enumerate}
\end{definition}

We are now finally in a position to state our Feynman-Kac formula for
the trace expectation, which is a special case of \cite[Theorem 4.9]{RigidityMSAO}:

\begin{lemma}
\label{Lemma: Raw Feynman-Kac}
For every $t>0$, $x\geq0$, and $1\leq i\leq r$, define
\begin{align}
\label{Equation: Trace Moment Off-Diagonal Term}
\mf M_t(\hat p,\hat U^i_t)=
\begin{cases}
1&\text{if }\hat N(t)=0,\\
\upsilon^{\hat N(t)}\mf C_t(\hat p,\hat U^i_t)&\text{otherwise},
\end{cases}
\end{align}
and
\begin{multline}
\label{Equation: Trace Moment Exponential Term}
\mf H_t(\hat U^i_t,X^{x,x}_t)=-\ka\int_0^{t}X^{x,x}_t(s)\d s\\
+\frac{(r-1)^2}2\|L_t(X^{x,x}_t)\|_2^2+\frac{\si^2}2\|L_t(\hat U^i_t,X^{x,x}_t)\|_\mu^2
-\sum_{j=1}^rw_j\mf L^{(j,0)}_{t}(\hat U^i_t,X^{x,x}_t).
\end{multline}
For every $t>0$, one has
\begin{align}
\label{Equation: Trace Moment Formula}
\mbf E\left[\mr{Tr}\big[\mr e^{-t\hat H_{\theta,\eta}}\big]\right]
=\int_{\mc A}
\Pi_X(t;x,x)\mbf E\left[\mf M_{t}(\hat p,\hat U^i_t)\mr e^{\mf H_{t}(\hat U^i_t,X^{x,x}_{t})}
\mbf 1_{\{\hat U^i_t(t)=i\}}\right]\d\mu(i,x).
\end{align}
\end{lemma}

\begin{remark}
Theorem 4.9 in \cite{RigidityMSAO} provides a formula for any joint moment
of the form $\mbf E\big[\mr{Tr}\big[\mr e^{-t_1\hat H}\big]\cdots\mr{Tr}\big[\mr e^{-t_n\hat H}\big]\big]$
for general vector-valued Schr\"odinger operators $\hat H=-\frac12\De+V+\xi$ on any one-dimensional
domain $I\subset\mbb R$. The statement in Lemma \ref{Lemma: Raw Feynman-Kac} can
be recovered from the latter by taking $n=1$ and $t=t_1$, then considering what is called "Case 2"
in \cite{RigidityMSAO} (i.e., an operator on the half-line), and then taking the potential function $V(i,x)=\ka x$.
\end{remark}

\begin{remark}
\label{Remark: Infinity times boundary means not touching}
The case where $w_j=\infty$
leads to some potential ambiguity in the definition of $\mf H_t(\hat U^i_t,X^{x,x}_t)$ in \eqref{Equation: Trace Moment Exponential Term}. Using the convention $\mr e^{-\infty\cdot c}=\mbf 1_{\{c=0\}}$ for $c\geq0$, we interpret the contribution
of $\mr e^{-\infty\cdot\mf L^{(j,0)}_{t}(\hat U^i_t,X^{x,x}_t)}$ in \eqref{Equation: Trace Moment Formula} as follows:
\[\mr e^{-\infty\cdot\mf L^{(j,0)}_{t}(\hat U^i_t,X^{x,x}_t)}
=\mbf 1_{\{\mf L^{(j,0)}_{t}(\hat U^i_t,X^{x,x}_t)=0\}}
=\mbf 1_{\{X^{x,x}_t(s)\neq0\text{ whenever }\hat U^i_t(s)=j\}}.\]
\end{remark}

\subsection{Proof of \eqref{Equation: Main Technical} Step 2. Simplifications}
\label{Section: Main Proof 1}

While \eqref{Equation: Trace Moment Formula} provides an exact formula
for the expected trace, it is quite difficult to parse due to its complexity.
Thus, our goal in this section
is to reformulate \eqref{Equation: Trace Moment Formula} in a way that is
amenable to computation, at least for the purposes of extracting the first- and second-order asymptotics.

For this purpose, we observe that most
of the complexity in \eqref{Equation: Trace Moment Formula} comes from the contribution of the path $\hat U^i_t$
to the combined local times \eqref{Equation: Combined Local Time 1}
and \eqref{Equation: Combined Local Time 2} and the combinatorics of $\mf C_t(\hat p,\hat U^i_t)$.
However, these terms are only difficult to understand when $\hat U^i_t$'s path has numerous jumps.
In this context,
our ability to calculate
the asymptotic \eqref{Equation: Main Technical} despite the complexity of
\eqref{Equation: Trace Moment Formula}
comes from the
fact that the leading order terms in the latter
only involve realizations of $\hat U^i_t$ that have no jump or two jumps, thus simplifying the analysis.

In fact, the unifying theme of the simplifications performed
in this section (as evidenced by Proposition \ref{Proposition: Reflected Bridge Asymptotics} below) is that
the exact calculation of leading terms and the estimate of remainder
terms in \eqref{Equation: Trace Moment Formula} can both be reduced to the small-time asymptotics of
three well-understood functionals of $X^{x,x}_t$ (neither or which involve $\hat U^{i}_t$),
namely, $\int_0^tX^{x,x}_t(s)\d s$, $\|L_t(X^{x,x}_t)\|_2$, and $\mf L_t^0(X^{x,x}_t)$. We postpone
the proofs of these asymptotics to Section \ref{Section: Bridge Asymptotics}.

To begin making this precise, we note that
we can use \eqref{Equation: Trace Moment Formula} to split
\begin{align}
\label{Equation: Trace into three terms}
\mbf E\left[\mr{Tr}\big[\mr e^{-t\hat H_{\theta,\eta}}\big]\right]=T_0(t)+T_2(t)+T_{\geq4}(t)
\end{align}
into a sum of expectations with
an indicator restricting the value of $\hat N(t)$ as follows:
\begin{align}
\label{Equation: N(t)=0 Case}
T_0(t)&=\int_{\mc A}\Pi_X(t;x,x)\mbf E\left[\mf M_{t}(\hat p,\hat U^i_t)\mr e^{\mf H_{t}(\hat U^i_t,X^{x,x}_{t})}
\mbf 1_{\{\hat U^i_t(t)=i\}\cap\{\hat N(t)=0\}}\right]\d\mu(i,x),\\
\label{Equation: N(t)=2 Case}
T_2(t)&=\int_{\mc A}\Pi_X(t;x,x)\mbf E\left[\mf M_{t}(\hat p,\hat U^i_t)\mr e^{\mf H_{t}(\hat U^i_t,X^{x,x}_{t})}
\mbf 1_{\{\hat U^i_t(t)=i\}\cap\{\hat N(t)=2\}}\right]\d\mu(i,x),\\
\label{Equation: N(t)>2 Case}
T_{\geq4}(t)&=\int_{\mc A}\Pi_X(t;x,x)\mbf E\left[\mf M_{t}(\hat p,\hat U^i_t)\mr e^{\mf H_{t}(\hat U^i_t,X^{x,x}_{t})}
\mbf 1_{\{\hat U^i_t(t)=i\}\cap\{\hat N(t)\geq4\}}\right]\d\mu(i,x).
\end{align}
We discuss the contribution
of each of these three terms to the $t\to0^+$ asymptotic in Sections \ref{Section: T0}--\ref{Section: T4},
and then summarize our findings in Section \ref{Section: Tk Summary}.

\subsubsection{The Constant Path}
\label{Section: T0}

We begin by simplifying $T_0(t)$ in \eqref{Equation: N(t)=0 Case}. This term is
the easiest to deal with, since if $\hat N(t)=0$, then $\hat U^i_t=i$ is constant on $[0,t]$.
Therefore, we get the following three immediate simplifications when $\hat N(t)=0$:
\begin{enumerate}[(1)]
\item $\mf M_t(\hat p,\hat U^i_t)=1$ by definition in \eqref{Equation: Trace Moment Off-Diagonal Term}.
\item Recalling \eqref{Equation: Combined Local Time 1} and \eqref{Equation: Combined Local Time 2},
we see that
\[\mf L^{(j,0)}_t(\hat U^i_t,X^{x,x}_t)=\begin{cases}
\mf L^0_t(X^{x,x}_t)&\text{if }j=i,\\
0&\text{if }j\neq i,
\end{cases}
\quad
and
\quad
L^{(j,y)}_t(\hat U^i_t,X^{x,x}_t)=\begin{cases}
L^y_t(X^{x,x}_t)&\text{if }j=i,\\
0&\text{if }j\neq i;
\end{cases}\]
hence \eqref{Equation: Trace Moment Exponential Term} simplifies to
\[
\mf H_t(\hat U^i_t,X^{x,x}_t)=
-\ka\int_0^{t}X^{x,x}_t(s)\d s
+\frac{(r-1)^2+\si^2}2\|L_t(X^{x,x}_t)\|_2^2
-w_i\mf L^0_{t}(X^{x,x}_t).\]
\item $\mbf 1_{\{\hat U^i_t(t)=i\}}=1$.
\end{enumerate}
If we combine all of these simplifications into
\eqref{Equation: N(t)=0 Case}, then we get
\begin{align}
\label{Equation: N(t)=0 Case 1}
T_0(t)
=\sum_{i=1}^r\int_0^\infty\Pi_X(t;x,x)\mbf E\left[\mr e^{-\ka\int_0^tX^{x,x}_t(s)\d s+\frac{(r-1)^2+\si^2}2\|L_t(X^{x,x}_t)\|_2^2-w_i\mf L^0_t(X^{x,x}_t)}\mbf 1_{\{\hat N(t)=0\}}\right]\d x.
\end{align}

By Definition \ref{Definition: Singular Process} (more specifically, the definition of $\hat N(t)$ in \eqref{Equation: Hat N Probability} therein),
\[\mbf E\big[\mbf 1_{\{\hat N(t)=0\}}\big|X^{x,x}_t\big]=\mbf P\big[\hat N(t)=0\big|X^{x,x}_t\big]=\mr e^{-\frac{(r-1)^2}2\|L_{t}(X^{x,x}_t)\|_2^2}.\]
Thus, by the tower property (conditioning on $X^{x,x}_t$),
\eqref{Equation: N(t)=0 Case 1} simplifies to
\begin{align}
\label{Equation: N(t)=0 Case 2}
T_0(t)
=\sum_{i=1}^r\int_0^\infty\Pi_X(t;x,x)\mbf E\left[\mr e^{-\ka\int_0^tX^{x,x}_t(s)\d s+\frac{\si^2}2\|L_t(X^{x,x}_t)\|_2^2-w_i\mf L^0_t(X^{x,x}_t)}\right]\d x.
\end{align}
In other words, if for $\msf w\in(-\infty,\infty]$ we denote
\begin{align}
\label{Equation: T0wi}
T_0^{(\msf w)}(t)
=\int_0^\infty\Pi_X(t;x,x)\mbf E\left[\mr e^{-\ka\int_0^tX^{x,x}_t(s)\d s+\frac{\si^2}2\|L_t(X^{x,x}_t)\|_2^2-\msf w\mf L^0_t(X^{x,x}_t)}\right]\d x,
\end{align}
then we have that
\begin{align}
\label{Equation: Split T0 into T0wi}
T_0(t)=T_0^{(w_1)}(t)+\cdots+T_0^{(w_r)}(t).
\end{align}
We now focus on understanding the asymptotics of each
$T_0^{(w_i)}(t)$ individually.

Toward this end, we first note that we expect
\begin{align}
\label{Equation: Informal Integral Approximation}
-\ka\int_0^tX^{x,x}_t(s)\d s\approx-\ka t x\qquad\text{as $t\to0^+$ for fixed }x>0.
\end{align}
Using this as a guide, we write
\[T_0^{(\msf w)}(t)
=\int_0^\infty\Pi_X(t;x,x)\mr e^{-\ka tx}\mbf E\left[\mr e^{-\ka\int_0^t(X^{x,x}_t(s)-x)\d s+\frac{\si^2}2\|L_t(X^{x,x}_t)\|_2^2-\msf w\mf L^0_t(X^{x,x}_t)}\right]\d x.\]
Next, by a Taylor expansion, we can write
\begin{align}
\label{Equation: e^z Order 2 Expansion}
\mr e^z=1+z+\msf R_2(z),\qquad\text{where }|\msf R_2(z)|\leq z^2\mr e^{|z|}.
\end{align}
Thus, we obtain
\begin{multline}
\label{Equation: Tw for real w}
T^{(\msf w)}_0(t)
=\int_0^\infty\Pi_X(t;x,x)\mr e^{-\ka tx}
\mbf E\Bigg[1
+\tfrac{\si^2}2\|L_t(X^{x,x}_t)\|_2^2\\
-\ka{\textstyle\int_0^t}\big(X^{x,x}_t(s)-x\big)\d s-\msf w\mf L^0_t(X^{x,x}_t)\\
+\msf R_2\Big(-\ka{\textstyle\int_0^t}\big(X^{x,x}_t(s)-x\big)\d s+\tfrac{\si^2}2\|L_t(X^{x,x}_t)\|_2^2-\msf w\mf L^0_t(X^{x,x}_t)\Big)\Bigg]\d x
\end{multline}
when $\msf w\in\mbb R$, and, using Remark \ref{Remark: Infinity times boundary means not touching}
to interpret $\mr e^{-\infty\cdot\mf L^0_t(X^{x,x}_t)}=\mbf 1_{\{\mf L^0_t(X^{x,x}_t)=0\}}$, we get
\begin{multline}
\label{Equation: Tinfinity}
T^{(\infty)}_0(t)
=\int_0^\infty\Pi_X(t;x,x)\mr e^{-\ka tx}
\mbf E\Bigg[\left(1
+\tfrac{\si^2}2\|L_t(X^{x,x}_t)\|_2^2\right)\mbf 1_{\{\mf L^0_t(X^{x,x}_t)=0\}}\\
-\ka\left({\textstyle\int_0^t}\big(X^{x,x}_t(s)-x\big)\d s\right)\mbf 1_{\{\mf L^0_t(X^{x,x}_t)=0\}}\\
+\msf R_2\Big(-\ka{\textstyle\int_0^t}\big(X^{x,x}_t(s)-x\big)\d s+\tfrac{\si^2}2\|L_t(X^{x,x}_t)\|_2^2\Big)\mbf 1_{\{\mf L^0_t(X^{x,x}_t)=0\}}\Bigg]\d x.
\end{multline}
With this in hand, we now state a number of asymptotics of the
functionals of $X^{x,x}_t$ appearing in \eqref{Equation: Tw for real w} and \eqref{Equation: Tinfinity},
whose proofs we postpone to Section \ref{Section: Bridge Asymptotics} below:

\begin{proposition}
\label{Proposition: Reflected Bridge Asymptotics}
As $t\to0^+$, we have the leading order contributions:
\begin{align}
\label{Equation: N(t)=0 Case - Leading Order}
&\int_0^\infty\Pi_X(t;x,x)\mr e^{-\ka tx}\d x
=\tfrac{1}{2\pi\ka}t^{-3/2}+\tfrac14+o(1),\\
\label{Equation: Self-Intersection}
&\int_0^\infty\Pi_X(t;x,x)\mr e^{-\ka tx}\mbf E\big[\|L_t(X^{x,x}_t)\|_2^2\big]\d x
=\tfrac{1}{2\ka}+o(1)
\end{align}
for the case $\msf w\in\mbb R$, and
\begin{align}
\label{Equation: N(t)=0 Case - Leading Order Dirichlet}
&\int_0^\infty\Pi_X(t;x,x)\mr e^{-\ka tx}\mbf E\big[\mbf 1_{\{\mf L^0_t(X^{x,x}_t)=0\}}\big]\d x
=\tfrac{1}{2\pi\ka}t^{-3/2}-\tfrac14+o(1),\\
\label{Equation: Self-Intersection Dirichlet}
&\int_0^\infty\Pi_X(t;x,x)\mr e^{-\ka tx}\mbf E\big[\|L_t(X^{x,x}_t)\|_2^2\mbf 1_{\{\mf L^0_t(X^{x,x}_t)=0\}}\big]\d x=\tfrac{1}{2\ka}+o(1)
\end{align}
for $\msf w=\infty$.
Moreover, as $t\to0^+$, we have the remainder terms
\begin{align}
\label{Equation: N(t)=0 Case - Integral}
&\int_0^\infty\Pi_X(t;x,x)\mr e^{-\ka tx}\mbf E\Big[{\textstyle\int_0^t}\big(X^{x,x}_t(s)-x\big)\d s\Big]\d x=o(1),\\
\label{Equation: N(t)=0 Case - Boundary LC}
&\int_0^\infty\Pi_X(t;x,x)\mr e^{-\ka tx}\mbf E\big[\mf L^0_t(X^{x,x}_t)\big]\d x=o(1),
\end{align}
and for every $\mf c>0$,
\begin{multline}
\label{Equation: N(t)=0 Case - Remainder}
\int_0^\infty\Pi_X(t;x,x)\mr e^{-\ka tx}\mbf E\Bigg[\Big(\big|{\textstyle\int_0^t}\big(X^{x,x}_t(s)-x\big)\d s\big|+\|L_t(X^{x,x}_t)\|_2+\mf L^0_t(X^{x,x}_t)\Big)^2\\
\cdot\exp\Bigg\{\mf c\Big(\big|{\textstyle\int_0^t}\big(X^{x,x}_t(s)-x\big)\d s\big|+\|L_t(X^{x,x}_t)\|_2^2+\mf L^0_t(X^{x,x}_t)\Big)\Bigg\}\Bigg]\d x=o(1).
\end{multline}
\end{proposition}

On the one hand,
if we apply \eqref{Equation: N(t)=0 Case - Leading Order} and \eqref{Equation: Self-Intersection}
to the first line of \eqref{Equation: Tw for real w}, then apply \eqref{Equation: N(t)=0 Case - Integral}
and \eqref{Equation: N(t)=0 Case - Boundary LC} to the second line of \eqref{Equation: Tw for real w},
and finally apply \eqref{Equation: N(t)=0 Case - Remainder} to the third line of \eqref{Equation: Tw for real w}
(as \eqref{Equation: N(t)=0 Case - Remainder} matches the upper bound on $\msf R_2$ in
\eqref{Equation: e^z Order 2 Expansion}),
then we get
\[T_0^{(\msf w)}(t)=\tfrac{1}{2\pi\ka}t^{-3/2}+\tfrac14\left(\tfrac{\si^2}{\ka}+1\right)+o(1)\qquad\text{as }t\to0^+\]
when $\msf w\in\mbb R$.
On the other hand,
if we apply \eqref{Equation: N(t)=0 Case - Leading Order Dirichlet} and \eqref{Equation: Self-Intersection Dirichlet}
to the first line of \eqref{Equation: Tinfinity}, and then apply \eqref{Equation: N(t)=0 Case - Integral}
and \eqref{Equation: N(t)=0 Case - Remainder} to the second and third lines of \eqref{Equation: Tinfinity}
respectively (using the trivial bound $\mbf 1_{\{\mf L^0_t(X^{x,x}_t)=0\}}\leq1$ in both cases), then we get
\[T_0^{(\infty)}(t)=\tfrac{1}{2\pi\ka}t^{-3/2}+
\tfrac14\left(\tfrac{\si^2}{\ka}-1\right)+o(1)\qquad\text{as }t\to0^+.\]
Given that
$r_0$ corresponds to the number of components $w_i$
that are finite (hence $r-r_0$ are the number of infinite components), then we conclude from \eqref{Equation: Split T0 into T0wi} that
\begin{align}
\label{Equation: N(t)=0 Case Combined}
T_0(t)=\tfrac{r}{2\pi\ka}t^{-3/2}+\tfrac14\big(\tfrac{r\si^2}{\ka}+2r_0-r\big)+o(1)\qquad\text{as }t\to0^+.
\end{align}

\subsubsection{Two Jumps}

We now simplify $T_2(t)$.
We begin by arguing that if $\hat U^i_t(t)=i$ and $\hat N(t)=2$,
then $\mf C_t(\hat p,\hat U^i_t)=1$, and therefore it follows from \eqref{Equation: Trace Moment Off-Diagonal Term} that
\begin{align}
\label{Equation: Mt when N(t)=2}
\mf M_t(\hat p,\hat U^i_t)=\upsilon^2.
\end{align}
To see this, we note that when $\hat U^i_t(t)=i$ and $\hat N(t)=2$,
the matching $\hat p$ can only be $\big\{\{1,2\}\big\}$,
and the two jumps in $\hat U^i_t(t)$'s path must be of the form $J_1=(i,j)$ and $J_2=J_1^*=(j,i)$
for some $j\neq i$.
When $\be=1,2$, this immediately implies that $\mf C_t(\hat p,\hat U^i_t)=1$
by Definition \ref{Definition: Combinatorial Constant}-(a) and -(b).
When $\be=4$, Definition \ref{Definition: Combinatorial Constant}-(c)
states that we must additionally argue that $\mf D_t(\hat p,\hat U^i_t)=1$.
For this, we note that since $J_1=J_2^*$,
the only two binary paths in $\mc B^{0,0}_2$ that respect $J_1$
and $J_2$ are $(0,0,0)$ and $(0,1,0)$, and $\{1,2\}$ cannot be a flip.
We then get from \eqref{Equation: D constant for H} that
$\mf D_t(\hat p,\hat U^i_t)=2^{-\hat N(t)/2}\cdot 2=1$, as desired.

Next, we note that by
\eqref{Equation: Hat N Probability}, \eqref{Equation: M Probability}, and the independence
of $M^i_k$ and $\{\hat N(t),X^{x,x}_t\}$,
\begin{align*}
\mbf P\big[\big\{\hat U^i_t(t)=i\big\}\cap\big\{\hat N(t)=2\big\}~\big|~X^{x,x}_t\big]
&=\mbf P\big[\big\{M^i_2=i\big\}\cap\big\{\hat N(t)=2\big\}~\big|~X^{x,x}_t\big]\\
&=\mbf P\big[M^i_2=i\big]\mbf P\big[\hat N(t)=2~\big|~X^{x,x}_t\big]\\
&=\frac1{(r-1)}\cdot\frac{(r-1)^2\|L_{t}(X^{x,x}_t)\|^2_2}{2}\mr e^{-\frac{(r-1)^2}2\|L_{t}(X^{x,x}_t)\|_2^2}\\
&=\frac{(r-1)\|L_{t}(X^{x,x}_t)\|^2_2}{2}\mr e^{-\frac{(r-1)^2}2\|L_{t}(X^{x,x}_t)\|_2^2}.
\end{align*}
Thus, if we denote the two-jump random walk bridge
\[\tilde U^{i,i}_{t}=\Big(\hat U^i_t~\Big|~\hat U^i_t(t)=i,~\hat N(t)=2\Big),\]
then an application of the tower property (conditional on $X^{x,x}_t$)
in \eqref{Equation: N(t)=2 Case} together with \eqref{Equation: Mt when N(t)=2} yields
\begin{multline}
\label{Equation: N(t)=2 Case 1}
T_2(t)=\frac{(r-1)\upsilon^2}2\sum_{i=1}^r\int_0^\infty\Pi_X(t;x,x)\mbf E\Bigg[\|L_{t}(X^{x,x}_t)\|^2_2\\
\cdot\exp\left(-\ka\int_0^tX^{x,x}_t(s)\d s+\frac{\si^2}2\|L_t(\tilde U^{i,i}_{t},X^{x,x}_t)\|_\mu^2-\sum_{j=1}^rw_j\mf L^{(j,0)}_{t}(\tilde U^{i,i}_{t},X^{x,x}_t)\right)\Bigg]\d x.
\end{multline}
Without loss of generality (see Corollary \ref{Corollary: Cancellation of Subcritical Coordinates}),
we can assume that
\begin{align}
\label{Equation: Ordering of wi's}
w_1<\cdots<w_{r_0}<\infty=w_{r_0+1}=\cdots=w_r,
\end{align}
keeping in mind that all or none of the $w_i$ could be equal to $\infty$
if $r_0=r$ or $r_0=0$. If we combine this with
the heuristic \eqref{Equation: Informal Integral Approximation},
then we can rewrite \eqref{Equation: N(t)=2 Case 1} as
\begin{multline}
\label{Equation: N(t)=2 Case 2}
T_2(t)=\frac{(r-1)\upsilon^2}2\sum_{i=1}^r\int_0^\infty\Pi_X(t;x,x)\mr e^{-\ka t x}\mbf E\Bigg[\|L_{t}(X^{x,x}_t)\|^2_2\\
\cdot\exp\left(-\ka\int_0^t\big(X^{x,x}_t(s)-x\big)\d s+\frac{\si^2}2\|L_t(\tilde U^{i,i}_{t},X^{x,x}_t)\|_\mu^2-\sum_{j\leq r_0}w_j\mf L^{(j,0)}_t(\tilde U^{i,i}_{t},X^{x,x}_t)\right)\\
\cdot\mbf 1_{\{\mc L^{(j,0)}_t(\tilde U^{i,i}_{t},X^{x,x}_t)=0,~\forall j>r_0\}}\Bigg]\d x.
\end{multline}

At this point, if we expand the exponential in the second line of \eqref{Equation: N(t)=2 Case 2} using
\begin{align}
\label{Equation: e^z Order 1 Expansion}
\mr e^z=1+\msf R_1(z),\qquad\text{where }|\msf R_1(z)|\leq|z|\mr e^{|z|},
\end{align}
then we see that
\begin{align}
\label{Equation: N(t)=2 Case 3}
T_2(t)=\frac{(r-1)\upsilon^2}2\sum_{i=1}^r\int_0^\infty\Pi_X(t;x,x)\mr e^{-\ka t x}\mbf E\Big[\mf I^{(1)}_{i,t}(x)\Big(1+\mf I^{(2)}_{i,t}(x)\Big)\Big]\d x,
\end{align}
where we define
\begin{align*}
\mf I^{(1)}_{i,t}(x)&=\|L_{t}(X^{x,x}_t)\|^2_2\mbf 1_{\{\mf L^{(j,0)}_t(\tilde U^{i,i}_{t},X^{x,x}_t)=0,~\forall j>r_0\}},\\
\mf I^{(2)}_{i,t}(x)&=\msf R_1\left(-\ka\int_0^t\big(X^{x,x}_t(s)-x\big)\d s+\frac{\si^2}2\|L_t(\tilde U^{i,i}_{t},X^{x,x}_t)\|_\mu^2-\sum_{j\leq r_0}w_j\mf L^{(j,0)}_t(\tilde U^{i,i}_{t},X^{x,x}_t)\right).
\end{align*}
We now analyze the contribution of each of these summands to $T_2(t)$'s asymptotics:

We begin with the contribution coming from $\mf I^{(1)}_{i,t}(x)$. For this, we notice that
\[\mbf 1_{\{\mf L^0_t(X^{x,x}_t)=0\}}\leq\mbf 1_{\{\mf L^{(j,0)}_t(\tilde U^{i,i}_{t},X^{x,x}_t)=0,~\forall j>r_0\}}\leq1,\]
where the first inequality follows from \eqref{Equation: Combined Local Time Sum}.
If we combine this with both \eqref{Equation: Self-Intersection} and \eqref{Equation: Self-Intersection Dirichlet},
then we obtain the $t\to0^+$ asymptotic
\begin{align}
\label{Equation: N(t)=2 Case 3.1}
\frac{(r-1)\upsilon^2}2\sum_{i=1}^r\int_0^\infty\Pi_X(t;x,x)\mr e^{-\ka t x}\mbf E\Big[\mf I^{(1)}_{i,t}(x)\Big]\d x
=\frac{r(r-1)\upsilon^2}{4\ka}+o(1).
\end{align}

Next, we analyze the contribution of $\mf I^{(1)}_{i,t}(x)\mf I^{(2)}_{i,t}(x)$.
On the one hand, we note that
\begin{multline}
\label{Equation: Local Time L2 Bound for N(t)=2}
\|L_t(\tilde U^{i,i}_{t},X^{x,x}_t)\|_\mu^2=\sum_{j=1}^r\int_0^\infty\left(L^{(j,y)}_t(\tilde U^{i,i}_{t},X^{x,x}_t)\right)^2\d y\\
\leq\int_0^\infty\left(\sum_{j=1}^rL^{(j,y)}_t(\tilde U^{i,i}_{t},X^{x,x}_t)\right)^2\d y
=\|L_t(X^{x,x}_t)\|_2^2,
\end{multline}
where the inequality follows from $z_1^2+\cdots+z_n^2\leq(z_1+\cdots+z_n)^2$
whenever $z_i\geq0$, and the last equality follows from \eqref{Equation: Combined Local Time Sum}.
On the other hand, if
\begin{align}
\label{Equation: wmin}
w_{\min}=\min\{w_1,\ldots,w_r\},
\end{align}
then we get from \eqref{Equation: Combined Local Time Sum} that
\begin{align}
\label{Equation: Local Time Boundary Bound for N(t)=2}
-\sum_{j\leq r_0}w_j\mf L^{(j,0)}_{t}(\tilde U^{i,i}_{t},X^{x,x}_t)\leq-w_{\min}\sum_{j=1}^r\mf L^{(j,0)}_{t}(\tilde U^{i,i}_{t},X^{x,x}_t)
=-w_{\min}\mf L^0_t(X^{x,x}_t).
\end{align}
Combining this with the trivial bound $\mbf 1_E\leq 1$ for any event $E$, we get that
\begin{multline}
\label{Equation: N(t)=2 Case 3.2 prelim}
\left|\mf I^{(1)}_{i,t}(x)\mf I^{(2)}_{i,t}(x)\right|\\
\leq C\|L_{t}(X^{x,x}_t)\|^2_2\msf R_1\Bigg(\mf c\Big(\big|{\textstyle\int_0^t}\big(X^{x,x}_t(s)-x\big)\d s\big|+\|L_t(X^{x,x}_t)\|_2^2+\mf L^0_t(X^{x,x}_t)\Big)\Bigg)
\end{multline}
for some constants $C,\mf c>0$.
Given that
\[\|L_{t}(X^{x,x}_t)\|^2_2\leq\big|{\textstyle\int_0^t}\big(X^{x,x}_t(s)-x\big)\d s\big|+\|L_t(X^{x,x}_t)\|_2^2+\mf L^0_t(X^{x,x}_t),\]
if we combine \eqref{Equation: N(t)=2 Case 3.2 prelim} with the upper bound on $\msf R_1$ in \eqref{Equation: e^z Order 1 Expansion},
then we get
\begin{multline*}
\left|\mf I^{(1)}_{i,t}(x)\mf I^{(2)}_{i,t}(x)\right|
\leq C\mf c\Big(\big|{\textstyle\int_0^t}\big(X^{x,x}_t(s)-x\big)\d s\big|+\|L_t(X^{x,x}_t)\|_2+\mf L^0_t(X^{x,x}_t)\Big)^2\\
\cdot\exp\Bigg\{\mf c\Big(\big|{\textstyle\int_0^t}\big(X^{x,x}_t(s)-x\big)\d s\big|+\|L_t(X^{x,x}_t)\|_2^2+\mf L^0_t(X^{x,x}_t)\Big)\Bigg\}.
\end{multline*}
Thus, an application of \eqref{Equation: N(t)=0 Case - Remainder} yields
\begin{align}
\label{Equation: N(t)=2 Case 3.2}
\frac{(r-1)\upsilon^2}2\sum_{i=1}^r\int_0^\infty\Pi_X(t;x,x)\mr e^{-\ka t x}\mbf E\Big[\mf I^{(1)}_{i,t}(x)\mf I^{(2)}_{i,t}(x)\Big]\d x=o(1).
\end{align}

Summarizing the argument for $T_2(t)$, if we combine \eqref{Equation: N(t)=2 Case 3.1} and \eqref{Equation: N(t)=2 Case 3.2}
with \eqref{Equation: N(t)=2 Case 3}, then we finally get
\begin{align}
\label{Equation: N(t)=2 Case 4}
T_2(t)=\tfrac{r(r-1)\upsilon^2}{4\ka}+o(1)\qquad \text{as }t\to0^+.
\end{align}

\subsubsection{Four or More Jumps}
\label{Section: T4}

We now discuss $T_{\geq4}(t)$.
For this purpose, we begin with the observation that
the combinatorial constant introduced in Definition \ref{Definition: Combinatorial Constant}
has a very simple upper bound (see \cite[(4.2)]{RigidityMSAO}):
\[|\mf C_t(\hat p,\hat U^i_t)|\leq 1.\]
Using essentially the same argument as in
\eqref{Equation: Local Time L2 Bound for N(t)=2}
and \eqref{Equation: Local Time Boundary Bound for N(t)=2}, we get the following:
\begin{align}
\label{Equation: Local Time L2 Bound General}
\|L_t(\hat U^i_t,X^{x,x}_t)\|_\mu^2\leq\|L_t(X^{x,x}_t)\|_2^2,
\end{align}
\begin{align}
\label{Equation: Local Time Boundary Bound General}
-\sum_{j=1}^rw_j\mf L^{(j,0)}_{t}(\hat U^i_t,X^{x,x}_t)\leq-w_{\min}\mf L^0_t(X^{x,x}_t),
\end{align}
where we recall that $w_{\min}$ is defined in \eqref{Equation: wmin}.
Combining all of this into \eqref{Equation: N(t)>2 Case}, we get from Jensen's inequality
(taking the absolute value inside the $\dd x$ integral) and $|\mf M_t(\hat p,\hat U^i_t)|=
\upsilon^{\hat N(t)}|\mf C_t(\hat p,\hat U^i_t)|$ that
\begin{multline}
\label{Equation: T4 1}
|T_{\geq4}(t)|\leq
\int_{\mc A}\Pi_X(t;x,x)\mbf E\Bigg[\mr e^{-\ka\int_0^tX^{x,x}_t(s)\d s+\frac{(r-1)^2+\si^2}2\|L_t(X^{x,x}_t)\|_2^2-w_{\min}\mf L^{0}_{t}(X^{x,x}_t)}\\
\cdot\upsilon^{\hat N(t)}
\mbf 1_{\{\hat U^i_t(t)=i\}\cap\{\hat N(t)\geq4\}}\Bigg]\d\mu(i,x).
\end{multline}

By \eqref{Equation: M Probability} and the independence
of $M^i_k$ and $\{\hat N(t),X^{x,x}_t\}$, we note that
\begin{align*}
\mbf E\left[\upsilon^{\hat N(t)}
\mbf 1_{\{\hat U^i_t(t)=i\}\cap\{\hat N(t)\geq4\}}\bigg|X^{x,x}_t\right]
=\frac1{1-r}\mbf E\left[\upsilon^{\hat N(t)}
\mbf 1_{\{\hat N(t)\geq4\}}\bigg|X^{x,x}_t\right].
\end{align*}
Then, by \eqref{Equation: Hat N Probability},
\begin{align*}
&\frac1{r-1}\mbf E\left[\upsilon^{\hat N(t)}
\mbf 1_{\{\hat N(t)\geq4\}}\bigg|X^{x,x}_t\right]\\
&=\frac1{r-1}\sum_{m=2}^\infty\upsilon^{2m}\left(\frac{(r-1)^2\|L_{t}(X^{x,x}_t)\|^2_2}{2}\right)^m\frac{1}{m!}\mr e^{-\frac{(r-1)^2}2\|L_{t}(X^{x,x}_t)\|_2^2}\\
&=\frac{\mr e^{-\frac{(r-1)^2}2\|L_{t}(X^{x,x}_t)\|_2^2}}{r-1}\left(\mr e^{\frac{(r-1)^2\upsilon^2}2\|L_{t}(X^{x,x}_t)\|_2^2}-1-\frac{(r-1)^2\upsilon^2}2\|L_{t}(X^{x,x}_t)\|_2^2\right),
\end{align*}
where the last line follows from $\mr e^z$'s Taylor expansion.
Given that $\mr e^z-1-z\leq z^2\mr e^{|z|}$ when $z\geq0$, this simplifies to
\begin{align}
\mbf E\left[\upsilon^{\hat N(t)}
\mbf 1_{\{\hat U^i_t(t)=i\}\cap\{\hat N(t)\geq4\}}\bigg|X^{x,x}_t\right]
\leq\frac{(r-1)^3\upsilon^4\|L_{t}(X^{x,x}_t)\|_2^4}{4}\mr e^{\frac{(\upsilon^2-1)(r-1)^2}2\|L_{t}(X^{x,x}_t)\|_2^2}.
\end{align}
Therefore, an application of the tower property in the left-hand side of \eqref{Equation: T4 1} yields
\begin{multline}
\label{Equation: Tgeq4 Final}
|T_{\geq4}(t)|\leq
\frac{r(r-1)^3\upsilon^4}{4}\int_{\mc A}\Pi_X(t;x,x)\mbf E\Bigg[\|L_t(X^{x,x}_t)\|_2^4\\
\cdot\mr e^{-\ka\int_0^tX^{x,x}_t(s)\d s+\frac{\upsilon^2(r-1)^2+\si^2}2\|L_t(X^{x,x}_t)\|_2^2-w_{\min}\mf L^{0}_{t}(X^{x,x}_t)}\Bigg]\d\mu(i,x).
\end{multline}

On the one hand, we note that there exists a large enough constant $\mf c>0$ such that
\begin{multline*}
\left(-\ka{\textstyle{\int}_0^tX^{x,x}_t(s)\d s+\frac{\upsilon^2(r-1)^2+\si^2}2\|L_t(X^{x,x}_t)\|_2^2-w_{\min}\mf L^{0}_{t}(X^{x,x}_t)}\right)\\
\leq \mf c\Big(\big|{\textstyle\int_0^t}\big(X^{x,x}_t(s)-x\big)\d s\big|+\|L_t(X^{x,x}_t)\|_2^2+\mf L^0_t(X^{x,x}_t)\Big).
\end{multline*}
On the other hand,
\[\|L_t(X^{x,x}_t)\|_2^4\leq\Big(\big|{\textstyle\int_0^t}\big(X^{x,x}_t(s)-x\big)\d s\big|+\|L_t(X^{x,x}_t)\|_2^2+\mf L^0_t(X^{x,x}_t)\Big)^2.\]
If we plug both of these estimates into \eqref{Equation: Tgeq4 Final},
then we get from \eqref{Equation: N(t)=0 Case - Remainder} that
\begin{align}
\label{Equation: N(t)geq4 Case}
T_{\geq 4}(t)=o(1)\qquad\text{as }t\to0^+.
\end{align}

\subsubsection{Summary}
\label{Section: Tk Summary}

Before moving on, we briefly mention that the combination of
\eqref{Equation: N(t)=0 Case Combined}
\eqref{Equation: N(t)=2 Case 4}, and
\eqref{Equation: N(t)geq4 Case}
implies the asymptotic \eqref{Equation: Main Technical}. Consequently, Theorem \ref{Theorem: Main Technical} now follows from
Proposition \ref{Proposition: Reflected Bridge Asymptotics}.
We now take on this final task.

\subsection{Proof of \eqref{Equation: Main Technical} Step 3. Reflected Brownian Bridge Asymptotics}
\label{Section: Bridge Asymptotics}

We now conclude this paper by proving Proposition \ref{Proposition: Reflected Bridge Asymptotics}.

\subsubsection{Preliminary: Brownian Motion Coupling}

Reflected Brownian motion bridges are more difficult to work with than Brownian bridges.
Thus, we begin the proof of Proposition \ref{Proposition: Reflected Bridge Asymptotics} by
examining how one can write (or estimate) expectations of functions of
$X^{x,x}_t$ as expectations of functions $B^{x,x}_t$ and/or $B^{x,-x}_t$ instead.
For this purpose, we recall that
in Definition \ref{Definition: B and X}, we couple $X$ with a standard
Brownian motion $B$ via $X=|B|$. Under this coupling,
we note that $X^{x,x}_t$ and a conditioned version of $B^x$ are related as follows:
\[X^{x,x}_t(s)=\big(|B^x(s)|~\big|~B^x(t)\in\{-x,x\}\big),\qquad s\in[0,t].\]
Given that
\[\mbf P\big[B^x(t)=\pm x~\big|~B^x(t)\in\{-x,x\}\big]=\frac{\Pi_B(t;x,\pm x)}{\Pi_B(t;x,x)+\Pi_B(t;x,-x)}=\frac{\Pi_B(t;x,\pm x)}{\Pi_X(t;x,x)}\]
and $\Pi_B(t;x,x)=\frac1{\sqrt{2\pi t}}$, $\Pi_B(t;x,-x)=\frac{\mr e^{-2x^2/t}}{\sqrt{2\pi t}}$,
for any functional $F$, we can write
\begin{multline}
\label{Equation: X to B coupling general}
\int_0^\infty\Pi_X(t;x,x)\mr e^{-\ka tx}\mbf E\big[F\big(X^{x,x}_t\big)\big]\d x\\
=\int_0^\infty\frac{\mr e^{-\ka tx}}{\sqrt{2\pi t}}\mbf E\big[F\big(|B^{x,x}_t|\big)\big]\d x
+\int_0^\infty\frac{\mr e^{-\ka tx-2x^2/t}}{\sqrt{2\pi t}}\mbf E\big[F\big(|B^{x,-x}_t|\big)\big]\d x.
\end{multline}

If we replace $F(X^{x,x}_t)$ by $F(X^{x,x}_t)\mbf 1_{\{\mf L^0_t(X^{x,x}_t)=0\}}$, then \eqref{Equation: X to B coupling general}
can be further simplified. Indeed, given that $B^{x,-x}_t$ must obviously
reach zero as it travels from $x$ to $-x$, one has $\mbf 1_{\{\mf L^0(|B^{x,-x}_t|)=0\}}=0$.
Thus, in that case \eqref{Equation: X to B coupling general} simplifies to
\begin{multline}
\label{Equation: X to B coupling Dirichlet 1}
\int_0^\infty\Pi_X(t;x,x)\mr e^{-\ka tx}\mbf E\big[F\big(X^{x,x}_t\big)\mbf 1_{\{\mf L^0_t(X^{x,x}_t)=0\}}\big]\d x\\
=\int_0^\infty\frac{\mr e^{-\ka tx}}{\sqrt{2\pi t}}\mbf E\left[F\big(|B^{x,x}_t|\big)\mbf 1_{\{\mf L^0(|B^{x,x}_t|)=0\}}\right]\d x.
\end{multline}
Next, if we write
\[\mbf 1_{\{\mf L^0(|B^{x,x}_t|)=0\}}=\mbf 1_{\{\forall s\leq t,~B^{x,x}_t(s)\neq0\}}=1-\mbf 1_{\{\exists s\leq t:~B^{x,x}_t(s)=0\}},\]
then we can write the right-hand side of \eqref{Equation: X to B coupling Dirichlet 1} as
\begin{align}
\label{Equation: X to B coupling Dirichlet 2}
\int_0^\infty\frac{\mr e^{-\ka tx}}{\sqrt{2\pi t}}\mbf E\left[F\big(|B^{x,x}_t|\big)\right]\d x
-\int_0^\infty\frac{\mr e^{-\ka tx}}{\sqrt{2\pi t}}\mbf E\left[F\big(|B^{x,x}_t|\big)\mbf 1_{\{\exists s\leq t:~B^{x,x}_t(s)=0\}}\right]\d x.
\end{align}
In order to simplify the second term in \eqref{Equation: X to B coupling Dirichlet 2}, we combine the following
observations:
\begin{enumerate}[(1)]
\item By the strong Markov property and the symmetry of Brownian motion about $0$,
we notice that
\[\big(|B^{x,x}_t|~\big|~\exists s\leq t:~B^{x,x}_t(s)=0\big)\deq |B^{x,-x}_t|\]
(i.e., we can couple the two processes above by reflecting the path of $B^{x,x}_t$ after its first passage to zero, thus making it end at $-x$ instead of $x$).
\item By using the joint density of the running minimum and the current value of a Brownian motion
(e.g., \cite[Chapter III, Exercise 3.14]{RevuzYor}), we obtain that
\[\mbf P\big[\exists s\leq t:~B^{x,x}_t(s)=0\big]=\mr e^{-2x^2/t}.\]
\end{enumerate}
Thus, if we turn the second expectation in \eqref{Equation: X to B coupling Dirichlet 2}
into a conditional expectation given $\{\exists s\leq t:~B^{x,x}_t(s)=0\}$
(by multiplying and dividing by the probability of that event, i.e., $\mr e^{-2x^2/t}$), then we
obtain the following counterpart to \eqref{Equation: X to B coupling general}
\begin{multline}
\label{Equation: X to B coupling Dirichlet}
\int_0^\infty\Pi_X(t;x,x)\mr e^{-\ka tx}\mbf E\big[F\big(X^{x,x}_t\big)\mbf 1_{\{\mf L^0_t(X^{x,x}_t)=0\}}\big]\d x\\
=\int_0^\infty\frac{\mr e^{-\ka tx}}{\sqrt{2\pi t}}\mbf E\big[F\big(|B^{x,x}_t|\big)\big]\d x
-\int_0^\infty\frac{\mr e^{-\ka tx-2x^2/t}}{\sqrt{2\pi t}}\mbf E\big[F\big(|B^{x,-x}_t|\big)\big]\d x.
\end{multline}

Finally, while \eqref{Equation: X to B coupling general} and \eqref{Equation: X to B coupling Dirichlet}
are useful for exact calculations (and thus necessary to obtain the explicit leading order terms in our
asymptotics), simpler expressions can be used when we are only interested in upper estimates.
For this purpose, we invoke \cite[(5.9)]{GaudreauLamarreEJP}, which states that for any nonnegative
functional $F$, one has
\begin{align}
\label{Equation: Nonnegative RBM to BM Estimate}
\mbf E\big[F(X^{x,x}_t)\big]\leq 2\mbf E\big[F(|B^{x,x}_t|)\big].
\end{align}
We are now prepared to undertake the proofs of our technical results.

\subsubsection{Proof of \eqref{Equation: N(t)=0 Case - Leading Order} and \eqref{Equation: N(t)=0 Case - Leading Order Dirichlet}}

By \eqref{Equation: X to B coupling general} and \eqref{Equation: X to B coupling Dirichlet} with $F=1$,
\[\int_0^\infty\Pi_X(t;x,x)\mr e^{-\ka tx}\d x
=\int_0^\infty\frac{\mr e^{-\ka tx}}{\sqrt{2\pi t}}\d x
+\int_0^\infty\frac{\mr e^{-\ka tx-2x^2/t}}{\sqrt{2\pi t}}\d x,\]
and
\[\int_0^\infty\Pi_X(t;x,x)\mr e^{-\ka tx}\mbf E\big[\mbf 1_{\{\mf L^0_t(X^{x,x}_t)=0\}}\big]\d x
=\int_0^\infty\frac{\mr e^{-\ka tx}}{\sqrt{2\pi t}}\d x
-\int_0^\infty\frac{\mr e^{-\ka tx-2x^2/t}}{\sqrt{2\pi t}}\d x.\]
We then get \eqref{Equation: N(t)=0 Case - Leading Order} and \eqref{Equation: N(t)=0 Case - Leading Order Dirichlet}
as follows: On the one hand,
\[\int_0^\infty\frac{\mr e^{-\ka tx}}{\sqrt{2\pi t}}\d x=\frac{1}{\sqrt{2\pi}\ka}t^{-3/2}.\]
On the other hand, by the change of variables $x\mapsto t^{1/2}x$ and dominated convergence,
\[\int_0^\infty\frac{\mr e^{-\ka tx-2x^2/t}}{\sqrt{2\pi t}}\d x
=\int_0^\infty\frac{\mr e^{-\ka t^{3/2}x-2x^2}}{\sqrt{2\pi}}\d x
=\int_0^\infty\frac{\mr e^{-2x^2}}{\sqrt{2\pi}}\d x+o(1)=\frac14+o(1).\]

\subsubsection{Proof of \eqref{Equation: Self-Intersection} and \eqref{Equation: Self-Intersection Dirichlet}}

By \eqref{Equation: X to B coupling general} and \eqref{Equation: X to B coupling Dirichlet},
the proof of \eqref{Equation: Self-Intersection} and \eqref{Equation: Self-Intersection Dirichlet} reduces to two claims:
As $t\to0^+$, one has
\begin{align}
\label{Equation: Self-Intersection Main}
&\int_0^\infty\frac{\mr e^{-\ka tx}}{\sqrt{2\pi t}}\mbf E\big[\|L_t(|B^{x,x}_t|)\|_2^2\big]\d x
=\tfrac{1}{2\ka}+o(1),\\
\label{Equation: Self-Intersection Remainder}
&\int_0^\infty\frac{\mr e^{-\ka tx-2x^2/t}}{\sqrt{2\pi t}}\mbf E\big[\|L_t(|B^{x,-x}_t|)\|_2^2\big]\d x
=o(1).
\end{align}
For this purpose, it is useful to note that, by Brownian scaling
and the change of variables $y\mapsto t^{1/2}y$,
\begin{multline}
\label{Equation: SILT Scaling}
\|L_t(|B^{x,x}_t|)\|_2^2
=\int_0^\infty L^y_t(|B^{x,\pm x}_t|)^2\d y
=t\int_0^\infty L^{t^{-1/2}y}_1(|B^{t^{-1/2}x,\pm t^{-1/2}x}_1|)^2\d y\\
=t^{3/2}\int_0^\infty L^{y}_1(|B^{t^{-1/2}x,\pm t^{-1/2}x}_1|)^2\d y=t^{3/2}\|L_1(|B^{t^{-1/2}x,\pm t^{-1/2}x}|)\|_2^2.
\end{multline}

We begin with the proof of \eqref{Equation: Self-Intersection Main}.
An application of \eqref{Equation: SILT Scaling} yields
\[\int_0^\infty\frac{\mr e^{-\ka tx}}{\sqrt{2\pi t}}\mbf E\big[\|L_t(|B^{x,x}_t|)\|_2^2\big]\d x
=\int_0^\infty\frac{t\,\mr e^{-\ka tx}}{\sqrt{2\pi}}\mbf E\big[\|L_1(|B^{t^{-1/2}x,t^{-1/2}x}_1|)\|_2^2\big]\d x.\]
If we now apply the change of variables $x\mapsto t^{-1}x$, then this becomes
\[\int_0^\infty\frac{\mr e^{-\ka tx}}{\sqrt{2\pi t}}\mbf E\big[\|L_t(|B^{x,x}_t|)\|_2^2\big]\d x
=\int_0^\infty\frac{\mr e^{-\ka x}}{\sqrt{2\pi}}\mbf E\big[\|L_1(|B^{t^{-3/2}x,t^{-3/2}x}_1|)\|_2^2\big]\d x.\]
The proof of \eqref{Equation: Self-Intersection Main} can thus be broken down into two steps:
\begin{multline}
\label{Equation: Self-Intersection Main 1}
\lim_{t\to0^+}\int_0^\infty\frac{\mr e^{-\ka x}}{\sqrt{2\pi}}\mbf E\big[\|L_1(|B^{t^{-3/2}x,t^{-3/2}x}_1|)\|_2^2\big]\d x\\
=\int_0^\infty\frac{\mr e^{-\ka x}}{\sqrt{2\pi}}\mbf E\big[\|L_1(B^{0,0}_1)\|_2^2\big]\d x
=\frac{\mbf E\big[\|L_1(B^{0,0}_1)\|_2^2\big]}{\sqrt{2\pi}\ka}
\end{multline}
and
\begin{align}
\label{Equation: Self-Intersection Main 2}
\mbf E\big[\|L_1(B^{0,0}_1)\|_2^2\big]=\tfrac{\sqrt{\pi}}{\sqrt 2}.
\end{align}

For \eqref{Equation: Self-Intersection Main 1}, we first notice that for any $y>0$,
\[L_1^y(|B^{t^{-3/2}x,t^{-3/2}x}_1|)=L_1^y(B^{t^{-3/2}x,t^{-3/2}x}_1)+L_1^{-y}(B^{t^{-3/2}x,t^{-3/2}x}_1).\]
Moreover, if we use the coupling
\[B^{t^{-3/2}x,t^{-3/2}x}_1=t^{-3/2}x+B^{0,0}_1,\]
then this becomes
\begin{align}
\label{Equation: Local Time of Absolute Value}
L_1^y(|B^{t^{-3/2}x,t^{-3/2}x}_1|)=L_1^{y-t^{-3/2}x}(B^{0,0}_1)+L_1^{-y-t^{-3/2}x}(B^{0,0}_1).
\end{align}
Therefore,
\begin{align}
\label{Equation: Self-Intersection Main 1.1}
\|L_1(|B^{t^{-3/2}x,t^{-3/2}x}_1|)\|_2^2=\int_0^\infty \left(L_1^{y-t^{-3/2}x}(B^{0,0}_1)+L_1^{-y-t^{-3/2}x}(B^{0,0}_1)\right)^2\d y.
\end{align}
For every fixed $x>0$, there exists some finite random $t_0$ small enough so that $B^{0,0}_1$'s path
does not intersect $(-\infty,-t^{-3/2}x]$ for every $t<t_0$. Therefore,
\[\lim_{t\to0^+}\|L_1(|B^{t^{-3/2}x,t^{-3/2}x}_1|)\|_2^2
=\lim_{t\to0^+}\int_0^\infty L_1^{y-t^{-3/2}x}(B^{0,0}_1)^2\d y=\|L_1(B^{0,0}_1)\|_2^2\qquad\text{almost surely.}\]
Consequently, \eqref{Equation: Self-Intersection Main 1} follows from dominated convergence if we prove that
\[\sup_{x>0}\mbf E\left[\sup_{t>0}\|L_1(|B^{t^{-3/2}x,t^{-3/2}x}_1|)\|_2^2\right]<\infty.\]
This immediately follows by combining \eqref{Equation: Self-Intersection Main 2}
with an application of the inequality $(z_1+z_1)^2\leq2(z_1^2+z_2^2)$ to \eqref{Equation: Self-Intersection Main 1.1},
the latter of which yields
\[\|L_1(|B^{t^{-3/2}x,t^{-3/2}x}_1|)\|_2^2\leq2\int_0^\infty L_1^{y-t^{-3/2}x}(B^{0,0}_1)^2+L_1^{-y-t^{-3/2}x}(B^{0,0}_1)^2\d y
=2\|L_t(B^{0,0}_1)\|_2^2.\]

With this in hand, the only element left in the proof of \eqref{Equation: Self-Intersection Main} is
\eqref{Equation: Self-Intersection Main 2}.
By Tonelli's theorem, we can write
\[\mbf E\big[\|L_t(B^{0,0}_1)\|_2^2\big]=\int_{\mbb R}\mbf E\big[L^y_1(B^{0,0}_1)^2\big]\d x.\]
By \cite[(1)]{Pitman}, we have the joint density
\begin{align}
\label{Equation: Pitman Joint Density}
\mbf P\left[L_1^y(B^0)\in\dd\ell,B^0(1)\in\dd x\right]=\frac{(|y|+|x-y|+\ell)\mr e^{-(|y|+|x-y|+\ell)^2/2}}{\sqrt{2\pi}}
\end{align}
for all $\ell>0$ and $x,y\in\mbb R$.
If we evaluate this at $x=0$, further noting that $\mbf P\left[B^0(1)\in\dd 0\right]=\Pi_B(1;0,0)=\frac{1}{\sqrt{2\pi}},$
then we get the conditional density
\[\mbf P\left[L_1^y(B^{0,0}_1)\in\dd\ell\right]=(2|y|+\ell)\mr e^{-(2|y|+\ell)^2/2}.\]
Therefore, we obtain \eqref{Equation: Self-Intersection Main 2} through a sequence of two
exact calculations:
\[\mbf E\left[L^y_1(B^{0,0}_1)^2\right]=\int_0^\infty\ell^2\cdot(2|y|+\ell)\mr e^{-(2|y|+\ell)^2/2}\d\ell
=2 \mr e^{-2 y^2}-2 \sqrt{2 \pi } |y|\, \mr{erfc}(\sqrt{2} |y|),\]
and then
\[\int_{\mbb R} 2 \mr e^{-2 y^2}-2 \sqrt{2 \pi } |y|\, \mr{erfc}(\sqrt{2} |y|)\d y=\tfrac{\sqrt{\pi}}{\sqrt 2},\]
where we recall that $\mr{erfc}(z)=1-\frac{2}{\sqrt{\pi}}\int_0^z\mr e^{-y^2}\d y$ is the complementary error function.

Now that \eqref{Equation: Self-Intersection Main} is proved,
the proof of \eqref{Equation: Self-Intersection} and \eqref{Equation: Self-Intersection Dirichlet}
relies on establishing \eqref{Equation: Self-Intersection Remainder}.
By \eqref{Equation: SILT Scaling}, this is equivalent to showing
\[t\int_0^\infty\frac{\mr e^{-\ka tx-2x^2/t}}{\sqrt{2\pi}}\mbf E\big[\|L_1(|B^{t^{-1/2}x,-t^{-1/2}x}_1|)\|_2^2\big]\d x=o(1).\]
If we apply the change of variables $x\mapsto t^{1/2}x$, then this becomes
\[t^{3/2}\int_0^\infty\frac{\mr e^{-\ka t^{3/2}x-2x^2}}{\sqrt{2\pi}}\mbf E\big[\|L_1(|B^{x,-x}_1|)\|_2^2\big]\d x=o(1).\]
As $\mr e^{-\ka t^{3/2}x}\to1$ as $t\to0$, by dominated convergence, we get \eqref{Equation: Self-Intersection Main}
if we show that
\begin{align}
\label{Equation: Self-Intersection Remainder 1}
\int_0^\infty\frac{\mr e^{-2x^2}}{\sqrt{2\pi}}\mbf E\big[\|L_1(|B^{x,-x}_1|)\|_2^2\big]\d x<\infty.
\end{align}
Using a similar argument as in \eqref{Equation: Local Time of Absolute Value}, we get
\begin{multline*}
\|L_1(|B^{x,-x}_1|)\|_2^2
=\int_0^\infty\Big(L^y_1(B^{x,-x}_1)+L^{-y}_1(B^{x,-x}_1)\Big)^2\d y\\
\leq2\int_0^\infty L^y_1(B^{x,-x}_1)^2+L^{-y}_1(B^{x,-x}_1)^2\d y
=2\|L_t(B^{x,-x})\|_2^2.
\end{multline*}
If we additionally use the fact that self-intersection local time is invariant
with respect to shifts in space, we get $\|L_t(B^{x,-x})\|_2^2=\|L_t(B^{0,-2x})\|_2^2$.
Thus \eqref{Equation: Self-Intersection Remainder 1} reduces to
\begin{align}
\label{Equation: Self-Intersection Remainder 2}
\int_0^\infty\frac{\mr e^{-2x^2}}{\sqrt{2\pi}}\mbf E\big[\|L_1(B^{0,-2x}_1)\|_2^2\big]\d x<\infty.
\end{align}
If we combine \eqref{Equation: Pitman Joint Density} with $\mbf P[B^0(1)\in\mr d(-2x)]=\Pi_B(1;0,-2x)=\frac{\mr e^{-2x^2}}{\sqrt{2\pi}}$,
then we get that
\[\mbf P\left[L_1^y(B^{0,-2x}_1)\in\dd\ell\right]=(|y|+|-2x-y|+\ell)\mr e^{-(|y|+|-2x-y|+\ell)^2/2}\mr e^{2x^2}.\]
We then get \eqref{Equation: Self-Intersection Remainder 2} through a sequence of exact calculations:
Firstly,
\begin{multline}
\label{Equation: Self-Intersection Remainder 3}
\mbf E\big[L^y_1(B^{0,-2x}_1)^2\big]=\int_0^\infty\ell^2\cdot(|y|+|-2x-y|+\ell)\mr e^{-(|y|+|-2x-y|+\ell)^2/2}\mr e^{2x^2}\d\ell\\
=\mr e^{2 x^2} \left(\sqrt{2 \pi } (| 2 x+y| +| y| ) \left(\mr{erf}\left(\frac{| 2 x+y| +| y| }{\sqrt{2}}\right)-1\right)+2\mr e^{-(| 2 x+y| +| y| )^2/2}\right),
\end{multline}
where $\mr{erf}(z)=\frac2{\sqrt\pi}\int_0^z\mr e^{-y^2}\d y$ denotes the error function.
Secondly, integrating the second line of \eqref{Equation: Self-Intersection Remainder 3} with
respect to $y$ on all of $\mbb R$ yields
\begin{multline}
\label{Equation: Self-Intersection Remainder 4}
\mbf E\big[\|L_1(B^{0,-2x}_1)\|_2^2\big]
=4 \sqrt{2 \pi }\mr e^{2 x^2} x^2 \mr{erf}(\sqrt{2} x)\\
+\tfrac{\sqrt\pi }{\sqrt2}\mr e^{2 x^2}(4 x^2+1) \mr{erfc}(\sqrt{2} x)
+2(1-2 \sqrt{2 \pi }\mr e^{2 x^2} x) x.
\end{multline}
Thirdly, integrating the right-hand side of \eqref{Equation: Self-Intersection Remainder 4}
multiplied by $\frac{\mr e^{-2x^2}}{\sqrt{2\pi}}$ finally yields
\[\int_0^\infty\frac{\mr e^{-2x^2}}{\sqrt{2\pi}}\mbf E\big[\|L_1(B^{0,-2x}_1)\|_2^2\big]\d x=\tfrac{\sqrt 2}{3\sqrt\pi}<\infty.\]
This proves \eqref{Equation: Self-Intersection Remainder 2}, thus concluding the proof
of \eqref{Equation: Self-Intersection} and \eqref{Equation: Self-Intersection Dirichlet}.

\subsubsection{Proof of \eqref{Equation: N(t)=0 Case - Integral}}

By \eqref{Equation: X to B coupling general}, we can split \eqref{Equation: N(t)=0 Case - Integral} into two parts:
\begin{align}
\label{N(t)=0 Case Part 1 - Integral 1}
&\int_0^\infty\frac{\mr e^{-\ka tx}}{\sqrt{2\pi t}}\mbf E\left[\int_0^t|B^{x,x}_t(s)|-x\d s\right]\d x=o(1),\\
\label{N(t)=0 Case Part 1 - Integral 2}
&\int_0^\infty\frac{\mr e^{-\ka tx-2x^2/t}}{\sqrt{2\pi t}}\mbf E\left[\int_0^t|B^{x,-x}_t(s)|-x\d s\right]\d x=o(1).
\end{align}
We begin with \eqref{N(t)=0 Case Part 1 - Integral 1}.
If we couple $B^{x,x}_t=x+B^{0,0}_t$ and then distribute the $t^{-1/2}$ term in $\frac1{\sqrt{2\pi t}}$ inside the integral, then
we get that (RHS is short for right-hand side)
\[\text{RHS of }\eqref{N(t)=0 Case Part 1 - Integral 1}
=\int_0^\infty\frac{\mr e^{-\ka tx}}{\sqrt{2\pi}}\mbf E\left[\int_0^tt^{-1/2}|x+B^{0,0}_t(s)|-t^{-1/2}x\d s\right]\d x.\]
Then, by a Brownian scaling and the changes of variables $s\mapsto t^{-1}s$ and $x\mapsto t^{-1/2}x$,
\[\text{RHS of }\eqref{N(t)=0 Case Part 1 - Integral 1}
=t^{3/2}\int_0^\infty\frac{\mr e^{-\ka t^{3/2}x}}{\sqrt{2\pi}}\mbf E\left[\int_0^1|x+B^{0,0}_1(s)|-x\d s\right]\d x.\]
At this point, if we show that
\begin{align}
\label{N(t)=0 Case Part 1 - Integral 1.1}
\int_0^\infty\left|\mbf E\left[\int_0^1|x+B^{0,0}_1(s)|-x\d s\right]\right|\d x=
\int_0^\infty\left|\int_0^1\mbf E\big[|x+B^{0,0}_1(s)|-x\big]\d s\right|\d x<\infty,
\end{align}
then we get \eqref{N(t)=0 Case Part 1 - Integral 1} (actually, the stronger statement RHS of \eqref{N(t)=0 Case Part 1 - Integral 1} $=O(t^{3/2})$) by dominated convergence
since $\mr e^{-\ka t^{3/2}x}\uparrow 1$ as $t\to0^+$.
Toward this end, since $B^{0,0}_1(s)$ is Gaussian with mean zero and variance $s(1-s)$, a direct calculation yields
\[\mbf E\big[|x+B^{0,0}_1(s)|-x\big]=\sqrt{\frac{2s(1-s)}{\pi}}\mr e^{-x^2/2s(1-s)}-x\,\mr{erfc}\left(\frac{x}{\sqrt{2s(1-s)}}\right).\]
Thus, by Jensen's inequality (taking the absolute value inside the $\dd s$ integral),
Tonelli's theorem (interchanging the $\dd x$ and $\dd s$ integrals), and direct calculations,
\begin{multline}
\text{RHS of }\eqref{N(t)=0 Case Part 1 - Integral 1.1}
\leq\int_0^1\int_0^\infty \sqrt{\frac{2s(1-s)}{\pi}}\mr e^{-x^2/2s(1-s)}+x\,\mr{erfc}\left(\frac{x}{\sqrt{2s(1-s)}}\right)\d x\dd s\\
=\frac32\int_0^1s(1-s)\dd s=\frac14.
\end{multline}
Thus proves \eqref{N(t)=0 Case Part 1 - Integral 1.1}, and therefore also \eqref{N(t)=0 Case Part 1 - Integral 1}.

We now move on to \eqref{N(t)=0 Case Part 1 - Integral 2}.
For this, we use the coupling
\[B^{x,-x}_t(s)=(1-\tfrac st)x-\tfrac stx+B^{0,0}_t(s).\]
By the reverse triangle inequality and the triangle inequality, for any $x>0$ and $s\in[0,t]$,
one has
\begin{multline}
\big||B^{x,-x}_t(s)|-x\big|
=\big||(1-\tfrac st)x-\tfrac stx+B^{0,0}_t(s)|-x\big|\\
\leq\big|(1-\tfrac st)x-\tfrac stx+B^{0,0}_t(s)-x\big|
=|B^{0,0}_t(s)-\tfrac{2s}tx|\leq|B^{0,0}_t(s)|+\tfrac{2s}{t}x.
\end{multline}
Combining this with $\mr e^{-\ka t x}\leq 1$, we get the upper bounds
\begin{align*}
|\text{RHS of }\eqref{N(t)=0 Case Part 1 - Integral 2}|
&\leq\int_0^\infty\frac{\mr e^{-2x^2/t}}{\sqrt{2\pi t}}\mbf E\left[\int_0^t|B^{0,0}_t(s)|+\tfrac{2s}{t}x\d s\right]\d x\\
&\leq \int_0^\infty\frac{\mr e^{-2x^2/t}}{\sqrt{2\pi t}}\cdot t\left(\mbf E\left[\sup_{0\leq s\leq t}|B^{0,0}_t(s)|\right]+x\right)\d x\\
&=\frac t4\mbf E\left[\sup_{0\leq s\leq t}|B^{0,0}_t(s)|\right]+\frac{t^{3/2}}{4\sqrt{2\pi}}\\
&=\frac{t^{3/2}}4\left(\mbf E\left[\sup_{0\leq s\leq 1}|B^{0,0}_1(s)|\right]+\frac{1}{\sqrt{2\pi}}\right),
\end{align*}
where the last line follows by Brownian scaling.
This proves \eqref{N(t)=0 Case Part 1 - Integral 2}.

\subsubsection{Proof of \eqref{Equation: N(t)=0 Case - Boundary LC}}

By \eqref{Equation: X to B coupling general}, \eqref{Equation: N(t)=0 Case - Boundary LC} reduces to
\begin{align}
\label{N(t)=0 Case Part 1 - Boundary LC 1}
\int_0^\infty\mr e^{-\ka tx}\Pi_B(t;x,\pm x)\mbf E\big[\mf L^0_t(B^{x,\pm x}_t)\big]\d x=o(1).
\end{align}
(Note that we use $B^{x,\pm x}_t$ instead of $|B^{x,\pm x}_t|$ above because
both processes have the same local time at zero.)
By Brownian scaling,
\[\mbf E\big[\mf L^0_t(B^{x,\pm x}_t)\big]=t^{1/2}\mbf E\big[\mf L^0_1(B^{t^{-1/2}x,\pm t^{-1/2}x}_1)\big]\]
Thanks to \cite[(1)]{Pitman}, for every $x,y\in\mbb R$ and $\ell>0$, we have that
\[\mbf P[\mf L^0_1(B^x)\in\dd \ell,B^x(1)\in\dd y]
=\frac{(|x|+|y|+\ell)}{\sqrt{2\pi}}\mr e^{-(|x|+|y|+\ell)^2/2}.\]
If we now turn the above into a conditional density given $B^x(1)=y$, then we get
\[\mbf P\big[\mf L^0_1(B^{x,y}_1)\in\dd \ell\big]
=\Pi_B(1;x,y)^{-1}(|x|+|y|+\ell)\mr e^{-(|x|+|y|+\ell)^2/2}.\]
By Brownian scaling,
\[\Pi_B(t;x,y)=\frac{\mr e^{-(x-y)^2/2t}}{\sqrt{2\pi t}}=t^{-1/2}\frac{\mr e^{-(t^{-1/2}x-t^{-1/2}y)^2/2}}{\sqrt{2\pi}}=t^{-1/2}\Pi_B(1;t^{-1/2}x,t^{-1/2}y);\]
therefore, we obtain that
\[\Pi_B(t;x,y)t^{1/2}\mbf P\big[\mf L^0_1(B^{t^{-1/2}x,t^{-1/2}y}_1)\in\dd \ell\big]
=(|t^{-1/2}x|+|t^{-1/2}y|+\ell)\mr e^{-(|t^{-1/2}x|+|t^{-1/2}y|+\ell)^2/2}.\]
With $y=\pm x$ for $x>0$, this then simplifies to
\begin{align}
\label{Equation: Boundary Local Time Density}
\Pi_B(t;x,\pm x)t^{1/2}\mbf P\big[\mf L^0_1(B^{t^{-1/2}x,\pm t^{-1/2}x}_1)\in\dd \ell\big]=(2t^{-1/2}x+\ell)\mr e^{-(2t^{-1/2}x+\ell)^2/2}.
\end{align}
Thus, an exact calculation yields
\begin{multline}
\Pi_B(t;x,\pm x)t^{1/2}\mbf E\big[\mf L^0_1(B^{t^{-1/2}x,t^{-1/2}\pm x}_1)\big]\\
=\int_0^\infty\ell\cdot(2t^{-1/2}x+\ell)\mr e^{-(2t^{-1/2}x+\ell)^2/2}\d\ell
=\frac{\sqrt\pi }{\sqrt2} \mr{erfc}\left(\sqrt 2 t^{-1/2}x\right).
\end{multline}
Consequently, another exact calculation yields
\[\text{RHS of }\eqref{N(t)=0 Case Part 1 - Boundary LC 1}=\int_0^\infty\mr e^{-\ka t x} \frac{\sqrt\pi }{\sqrt2} \mr{erfc}\left(\sqrt 2 t^{-1/2}x\right) \d x=\frac{\sqrt{\frac{\pi }{2}} \left(1-e^{\ka^2 t^3/8} \text{erfc}\left(\frac{\ka t^{3/2}}{2 \sqrt{2}}\right)\right)}{\ka t}\]
We now conclude the proof of \eqref{N(t)=0 Case Part 1 - Boundary LC 1} thanks to the asymptotic
\[\frac{\sqrt{\frac{\pi }{2}} \left(1-e^{\ka^2 t^3/8} \text{erfc}\left(\frac{\ka t^{3/2}}{2 \sqrt{2}}\right)\right)}{\ka t}=O(t^{1/2})\qquad\text{as }t\to0^+,\]
which follows from the product of the Taylor expansions $\mr{erfc}(z)=1-\frac{2z}{\sqrt\pi}+O(z^2)$
and $\mr e^z=1+O(z)$ as $z\to0$.

\subsubsection{Proof of \eqref{Equation: N(t)=0 Case - Remainder}}

With \eqref{Equation: N(t)=0 Case - Leading Order}--\eqref{Equation: N(t)=0 Case - Boundary LC} proved,
the only estimate left to prove in the paper is \eqref{Equation: N(t)=0 Case - Remainder}.
If we combine $\Pi_X(t;x,x)\leq\frac{2}{\sqrt{2\pi t}}$ and \eqref{Equation: Nonnegative RBM to BM Estimate},
then \eqref{Equation: N(t)=0 Case - Remainder} reduces to
\begin{multline}
\label{Equation: Last Remainder Bridge 1}
t^{-1/2}\int_0^\infty\mr e^{-\ka tx}\mbf E\Bigg[\Big(\big|{\textstyle\int_0^t}\big(|B^{x,x}_t(s)|-x\big)\d s\big|+\|L_t(|B^{x,x}_t|)\|_2+\mf L^0_t(B^{x,x}_t)\Big)^2\\
\cdot\exp\Bigg\{\mf c\Big(\big|{\textstyle\int_0^t}\big(|B^{x,x}_t(s)|-x\big)\d s\big|+\|L_t(|B^{x,x}_t|)\|_2^2+\mf L^0_t(B^{x,x}_t)\Big)\Bigg\}\Bigg]\d x=o(1),
\end{multline}
where there is no absolute value on $B^{x,x}_t$ in the boundary local time $\mf L^0_t$, since $B^{x,x}_t$ and $|B^{x,x}_t|$ both
have the same local time at zero.
If we use the coupling $B^{x,x}_t=x+B^{0,0}_t$ and the reverse triangle inequality, then we get
\[\big|{\textstyle\int_0^t}\big(|B^{x,x}_t(s)|-x\big)\d s\big|\leq{\textstyle\int_0^t}\big||x+B^{0,0}_t(s)|-x\big|\d s\leq t\sup_{0\leq s\leq t}|B^{0,0}_t(s)|.\]
If we combine this with an application of H\"older's inequality, then \eqref{Equation: Last Remainder Bridge 1} reduces to
\begin{multline*}
t^{-1/2}\int_0^\infty\mr e^{-\ka tx}\mbf E\Bigg[\Big(t\sup_{0\leq s\leq t}|B^{0,0}_t(s)|+\|L_t(|B^{x,x}_t|)\|_2+\mf L^0_t(B^{x,x}_t)\Big)^4\Bigg]^{1/2}\\
\cdot\mbf E\Bigg[\exp\Bigg\{2\mf c\Big(t\sup_{0\leq s\leq t}|B^{0,0}_t(s)|+\|L_t(|B^{x,x}_t|)\|_2^2+\mf L^0_t(B^{x,x}_t)\Big)\Bigg\}\Bigg]^{1/2}\d x=o(1).
\end{multline*}
If we then apply a Brownian scaling, then we get the further reduction
\begin{multline}
\label{Equation: Last Remainder Bridge 2}
t^{-1/2}\int_0^\infty\mr e^{-\ka tx}\mbf E\Bigg[\Big(t^{3/2}\sup_{0\leq s\leq 1}|B^{0,0}_1(s)|+t^{3/2}\|L_1(|B^{x_t,x_t}_1|)\|_2+t^{1/2}\mf L^0_1(B^{x_t,x_t}_1)\Big)^4\Bigg]^{1/2}\\
\cdot\mbf E\Bigg[\exp\Bigg\{2\mf c\Big(t^{3/2}\sup_{0\leq s\leq 1}|B^{0,0}_1(s)|+t^{3/2}\|L_1(|B^{x_t,x_t}_1|)\|_2^2+t^{1/2}\mf L^0_1(B^{x_t,x_t}_1)\Big)\Bigg\}\Bigg]^{1/2}\d x=o(1),
\end{multline}
where we denote $x_t=t^{-1/2}x$ for simplicity.

By H\"older's inequality,
\begin{align}
\nonumber
\sup_{x>0,~t\in(0,1]}\mbf E\Bigg[\exp\Bigg\{2\mf c\Big(t^{3/2}\sup_{0\leq s\leq 1}|B^{0,0}_1(s)|+t^{3/2}\|L_1(|B^{x_t,x_t}_1|)\|_2^2+t^{1/2}\mf L^0_1(B^{x_t,x_t}_1)\Big)\Bigg\}\Bigg]^{1/2}&\\
\label{Equation: Last Remainder Bridge 3}
\leq
\mbf E\left[\mr e^{6\mf c\sup_{0\leq s\leq 1}|B^{0,0}_1(s)|}\right]^{1/6}
\sup_{z>0}\mbf E\left[\mr e^{6\mf c\|L_1(|B^{z,z}_1|)\|_2^2}\right]^{1/6}
\sup_{z>0}\mbf E\left[\mr e^{6\mf c\mf L^0_1(B^{z,z}_1)}\right]^{1/6}.&
\end{align}
Since $|B^{0,0}_1|$ is a Bessel bridge of dimension 1,
the first expectation in the second line of \eqref{Equation: Last Remainder Bridge 3} is finite thanks to the Gaussian
tail bound for Bessel bridge maxima in \cite[Remark 3.1]{GruetShi}. The third expectation in the second line of
\eqref{Equation: Last Remainder Bridge 3} is finite thanks to \cite[Lemma 5.8]{GaudreauLamarreEJP}.
Finally, if we use the bound $(z_1+z_2)^2\leq2(z_1^2+z_2^2)$, then we see that
\begin{multline*}
\|L_1(|B^{z,z}_1|)\|_2^2
=\int_0^\infty L^y_1(|B^{z,z}_1|)^2\d y
=\int_0^\infty\Big(L^y_1(B^{z,z}_1)+L^{-y}_1(B^{z,z}_1)\Big)^2\d y\\
\leq2\int_0^\infty L^y_1(B^{z,z}_1)^2+L^{-y}_1(B^{z,z}_1)^2\d y=2\|L_1(B^{z,z}_t)\|_2^2,
\end{multline*}
hence
\begin{align}
\label{Equation: SILT Exponential Moments}
\sup_{z>0}\mbf E\left[\mr e^{6\mf c\|L_1(|B^{z,z}_1|)\|_2^2}\right]^{1/6}\leq\sup_{z>0}\mbf E\left[\mr e^{12\mf c\|L_1(B^{z,z}_1)\|_2^2}\right]^{1/6}<\infty,
\end{align}
where the last inequality is from \cite[Lemma 5.11]{GaudreauLamarreEJP}.
In summary, \eqref{Equation: Last Remainder Bridge 3} is finite, which means that we can now reduce \eqref{Equation: Last Remainder Bridge 2} to
\begin{align}
\label{Equation: Last Remainder Bridge 4}
t^{-1/2}\int_0^\infty\mr e^{-\ka tx}\mbf E\Bigg[\Big(t^{3/2}\sup_{0\leq s\leq 1}|B^{0,0}_1(s)|+t^{3/2}\|L_1(|B^{x_t,x_t}_1|)\|_2^2+t^{1/2}\mf L^0_1(B^{x_t,x_t}_1)\Big)^4\Bigg]^{1/2}\d x=o(1).
\end{align}

Since there exists some $c>0$ such that $(z_1+z_2+z_3)^4\leq c(z_1^4+z_2^4+z_3^4)$ for all $z_i\geq0$ (by Jensen's inequality),
and $\sqrt{z_1+z_2+z_3}\leq\sqrt z_1+\sqrt z_2+\sqrt z_3$ for all $z_i\geq0$, we can further reduce \eqref{Equation: Last Remainder Bridge 4} to
\begin{multline}
\label{Equation: Last Remainder Bridge 5}
t^{-1/2}\int_0^\infty\mr e^{-\ka tx}\Bigg(\mbf E\Big[t^6\sup_{0\leq s\leq 1}|B^{0,0}_1(s)|^4\Big]^{1/2}+\mbf E\Big[t^6\|L_1(|B^{x_t,x_t}_1|)\|_2^8\Big]^{1/2}\\
+\mbf E\big[t^2\mf L^0_1(B^{x_t,x_t}_1)^4\big]^{1/2}\Bigg)=o(1).
\end{multline}
The first line of \eqref{Equation: Last Remainder Bridge 5} vanishes thanks to the combination of three facts: Firstly,
\[t^{-1/2}\int_0^\infty\mr e^{-\ka t x}(t^{6})^{1/2}\d x=\frac{t^{3/2}}{\ka}=o(1).\]
Secondly, 
\[\mbf E\Big[\sup_{0\leq s\leq 1}|B^{0,0}_1(s)|^4\Big]<\infty\]
thanks to the fact that $|B^{0,0}_1|$ supremum has Gaussian tails (\cite[Remark 3.1]{GruetShi}).
Thirdly, it follows from \eqref{Equation: SILT Exponential Moments} that
\[\sup_{z>0}\mbf E\big[\|L_1(|B^{z,z}_1|)\|_2^8\big]<\infty.\]
Thus, we now only need to show that the second line of \eqref{Equation: Last Remainder Bridge 5} vanishes, which we simplify as
as follows:
\begin{align}
\label{Equation: Last Remainder Bridge 6}
t^{1/2}\int_0^\infty\mr e^{-\ka tx}\mbf E\big[\mf L^0_1(B^{t^{-1/2}x,t^{-1/2}x}_1)^4\big]^{1/2}\d x=o(1).
\end{align}

Toward this end, noting that $\ka t\mr e^{-\ka t x}\dd x$ is a probability measure on $[0,\infty)$ and that the square root function is concave, it follows from Jensen's inequality that
\[t^{1/2}\int_0^\infty\mr e^{-\ka tx}\mbf E\big[\mf L^0_1(B^{t^{-1/2}x,t^{-1/2}x}_1)^4\big]^{1/2}\d x
\leq \frac{t^{1/2}}{\ka t}\left(\int_0^\infty\ka t\,\mr e^{-\ka tx}\mbf E\big[\mf L^0_1(B^{t^{-1/2}x,t^{-1/2}x}_1)^4\big]\d x\right)^{1/2}.\]
If we now group together the powers of $t$ in the above, we can reduce
\eqref{Equation: Last Remainder Bridge 6} to
\begin{align}
\label{Equation: Last Remainder Bridge 7}
\int_0^\infty\mr e^{-\ka tx}\mbf E\big[\mf L^0_1(B^{t^{-1/2}x,t^{-1/2}x}_1)^4\big]\d x=o(1).
\end{align}
If we combine \eqref{Equation: Boundary Local Time Density} with the fact that
$\Pi(t;x,x)t^{1/2}=\frac{1}{\sqrt{2\pi}}$, then we get from an exact calculation that
\begin{multline*}
\mbf E\big[\mf L^0_1(B^{t^{-1/2}x,t^{-1/2}x}_1)^4\big]
=\int_0^\infty \ell^4\cdot\sqrt{2\pi}(2t^{-1/2}x+\ell)\mr e^{-(2t^{-1/2}x+\ell)^2/2}\d\ell\\
=8 \sqrt{2 \pi } e^{-2 x^2/t}(t+2 x^2)t^{-1}-8\pi(3 t x+4 x^3)\mr{erfc}\left(\tfrac{\sqrt{2} x}{\sqrt{t}}\right)t^{-3/2}.
\end{multline*}
If we integrate the second line above multiplied by $\mr e^{-\ka t x}$, an exact
calculation yields
\begin{multline}
\label{Equation: Last Remainder Bridge 8}
\int_0^\infty\mr e^{-\ka tx}\mbf E\big[\mf L^0_1(B^{t^{-1/2}x,t^{-1/2}x}_1)^4\big]\d x\\
=\frac{192 \pi\mr e^{\ka^2 t^3/8} \mr{erfc}\left(\frac{\ka t^{3/2}}{2 \sqrt{2}}\right)}{\ka^4 t^{11/2}}-\frac{192 \pi }{\ka^4 t^{11/2}}+\frac{96 \sqrt{2 \pi }}{\ka^3 t^4}-\frac{24 \pi }{\ka^2 t^{5/2}}+\frac{8 \sqrt{2 \pi }}{\ka t}.
\end{multline}
To  see why this vanishes,
we combine the Taylor expansions (as $z\to0$)
\[\mr e^{z}=1+z+O(z^2)
\qquad\text{and}\qquad
\mr{erfc}(z)=1-\tfrac{2}{\sqrt{\pi}}z+\tfrac{2}{3 \sqrt{\pi }}z^3+O(z^5),\]
which yields the $t\to0^+$ asymptotics
\[\mr e^{\ka^2 t^3/8}=1+\tfrac{\ka^2 t^3}{8}+O(t^6)
\qquad\text{and}\qquad
\mr{erfc}\left(\frac{\ka t^{3/2}}{2 \sqrt{2}}\right)=1-\tfrac{\ka t^{3/2}}{\sqrt{2 \pi }}+\tfrac{\ka^3 t^{9/2}}{24 \sqrt{2 \pi }}+O(t^{15/2}).\]
If we take a product of the above, and then multiply the result by $\frac{192 \pi}{\ka^4 t^{11/2}}
$, then we get
\[\frac{192 \pi\mr e^{\ka^2 t^3/8} \mr{erfc}\left(\frac{\ka t^{3/2}}{2 \sqrt{2}}\right)}{\ka^4 t^{11/2}}
=\frac{192 \pi }{\ka^4 t^{11/2}}-\frac{96 \sqrt{2 \pi }}{\ka^3 t^4}+\frac{24 \pi }{\ka^2 t^{5/2}}-\frac{8 \sqrt{2 \pi }}{\ka t}+O(t^{1/2}).\]
If we put this into \eqref{Equation: Last Remainder Bridge 8},
then we obtain \eqref{Equation: Last Remainder Bridge 6}.
This finally concludes the proof of \eqref{Equation: N(t)=0 Case - Remainder},
and thus of Proposition \ref{Proposition: Reflected Bridge Asymptotics}.
Consequently, Theorem \ref{Theorem: Main Technical} is now proved.

\appendix
\section{Heuristic \ref{Heuristic: beta-ensembles}}
\label{Section: beta-ensembles heuristic}

Recall the definitions of $Z_{\be,n}$ and $H_n(\mc B^{\be,n})$ in \eqref{Equation: Beta Ensemble Density}
and \eqref{Equation: Beta Ensembles Energy}. Straightforward calculations using
\eqref{Equation: Beta Ensemble Density} reveal that
\begin{align}
\label{Equation: Energy to Partition Function}
\mbf E\big[H_n(\mc B^{\be,n})\big]=-\tfrac{1}{n}\tfrac{\partial}{\partial\be}\log Z_{\be,n}
\qquad\text{and}\qquad
\mbf{Var}\big[H_n(\mc B^{\be,n})\big]=\tfrac{1}{n^2}\tfrac{\partial^2}{\partial^2\be}\log Z_{\be,n}.
\end{align}
In the special case where $V(x)=x^2/2$, a classical calculation using a Selberg integral \cite{Selberg}
yields the explicit formula $Z_{\be,n}=Z_{G\be E,n}$, where we define
\[Z_{G\be E,n}=(2\pi)^{n/2}(n\be/2)^{-\be n^2/4+(\be/4-1/2)n}\frac{\prod_{j=1}^n\Gamma(1+j\be/2)}{\Gamma(1+\be/2)^n}.\]
Asymptotics using several applications of the Mathematica "Limit" function reveal that, as $n\to\infty$, one has
\[-\tfrac{1}{n}\tfrac{\partial}{\partial\be}\log Z_{G\be E,n}=\tfrac38n-\tfrac12\log n-\tfrac{1}{4}\left(1+\log(\be^2/4)-2\tfrac{\Gamma'(1+\be/2)}{\Gamma(1+\be/2)}\right)+o(1)\]
and $\tfrac{1}{n^2}\tfrac{\partial^2}{\partial^2\be}\log Z_{G\be E,n}=o(1)$. Combining this with \eqref{Equation: Energy to Partition Function},
we get that when $V(x)=x^2/2$, the limit \eqref{Equation: Energy Limit in Probability} holds.
Using the well-known identity between the polygamma functions and the Hurwitz zeta function (e.g., \cite[6.4.10]{AbramowitzStegun}), we have that
for $\be>0$,
\[\frac{\mr d}{\mr d\be}\left(\log(\be^2/4)-2\frac{\Gamma'(1+\be/2)}{\Gamma(1+\be/2)}\right)=\frac{2}{\be}-\sum_{k=0}^\infty\frac{1}{(k+1+\be/2)^2}
>\frac{2}{\be}-\int_{-1}^\infty\frac{1}{(k+1+\be/2)^2}\d k=0,\]
where the inequality follows from the integral test since $k\mapsto(k+1+\be/2)^{-2}$ is strictly decreasing.
This implies that $\be\mapsto\log(\be^2/4)-2\tfrac{\Gamma'(1+\be/2)}{\Gamma(1+\be/2)}$ is invertible for $\be\in(0,\infty)$.

\bibliographystyle{plain}
\bibliography{Bibliography}

\end{document}